\documentclass[11pt,a4paper,equation]{article}
\usepackage[greek,american]{babel}
\usepackage[iso-8859-7]{inputenc}
\usepackage{amssymb,amsfonts}
\usepackage[dvips]{graphicx}
\usepackage{amsmath}
\usepackage{epsfig}
\usepackage{latexsym}
\usepackage{makeidx}
\usepackage{pst-math,pst-xkey}

\numberwithin{equation}{section}
\makeindex
\begin{document}
\date{}
\author{Aristides V. Doumas\\
Department of Mathematics\\
National Technical University of Athens\\
Zografou Campus\\
157 80 Athens, GREECE\\
\underline{aris.doumas@hotmail.com}\\\\
and\\\\
Vassilis G. Papanicolaou\\
Department of Mathematics\\
National Technical University of Athens\\
Zografou Campus\\
157 80 Athens, GREECE\\
\underline{papanico@math.ntua.gr}}
\title{The Siblings of the Coupon Collector}
\maketitle
\begin{abstract}
The following variant of the collector's problem has attracted considerable
attention relatively recently (see, e.g., \cite{P}, \cite{FO}, \cite{FD},
\cite{Ad}, and \cite{Ro}, listed here in chronological order): There is one
main collector who collects coupons. Assume there are $N$ different types of
coupons with, in general, unequal occurring probabilities. When the main
collector gets a "double", she gives it to her older brother; when this
brother gets a "double", he gives it to the next brother, and so on. Hence,
when the main collector completes her collection, the album of the $j$-th
sibling, $j = 2, 3, \dots$, will still have $U_j^N$ empty spaces. In this article we develop techniques of computing asymptotics of the average
$E[U_j^N]$ of $U_j^N$ as $N\rightarrow \infty,$ for a large class of families
of coupon probabilities. We also give various illustrative examples.

\end{abstract}
\textbf{Keywords.} Urn problems; (generalized) coupon collector's problem (GCCP); asymptotics.\\\\
\textbf{2010 AMS Mathematics Classification.} 60F05; 60F99.
\section{Introduction}

\subsection{Preliminaries}

The classical \textbf{``coupon collector's problem'' (CCP)} concerns a population
(e.g. fishes, viruses, genes, words, baseball cards, etc.)
whose members are of $N$ different \emph{types}. The members of the population are
sampled independently with replacement and their types are recorded. CCP pertains to the family of urn problems along with other famous problems, such as the birthday, or occupancy. Its origin can be traced back to De Moivre's treatise \textit{De Mensura Sortis} of 1712 (see, e.g., \cite{Ho}) and Laplace's pioneering work \textit{Theorie Analytique de Probabilites} of 1812 (see \cite{D-H}). The problem became popular in the 1930's when the Dixie Cup company introduced a highly successful program by which children collected Dixie lids to receive ``Premiums,'' beginning with illustrations of their favored Dixie Circus characters, and then Hollywood stars and major league baseball players (for the Dixie Cup company history see \cite{DC}).\\
For $1\leq k\leq N$ let $p_{k}$ be the probability that a member of the population is of type $k$, $\sum_{k=1}^{N}p_{k}=1$. Let $T_{N}$ be the number of trials it takes until all $N$ types are detected (at least once).\\
General results for the simplest CCP (i.e. the case of equal probabilities) had appeared in some relatively unknown works (see \cite{M}, \cite{G}, where the entertaining term \textit{cartophily} appeared in the title of these papers). However, some of the  classical references for this case of the CCP are W. Feller's well known work \cite{F}, D.J. Newman's and L. Shepp's paper on the \textit{Double Dixie Cup} problem \cite{N-S} (where they answered the question: how long, on average does it take to obtain $m$ complete sets of $N$ coupons), and a paper of P. Erd\H{o}s and A. R$\acute{e}$nyi, where the limit distribution of the random variable $T_{N}$ has been given (see \cite{E}). Since then, CCP has attracted the attention of various researchers due to the fact that it has found many applications in many areas of science (computer science/search algorithms, mathematical programming, optimization, learning processes, engineering, ecology, as well as linguistics---see, e.g., \cite{B-H}, \cite{Ma}).\\
For the general case of unequal probabilities regarding the asymptotics of the moments, as well as for the limit distribution of $T_{N}$, there is a plethora of referenceses (see for instance \cite{B-B}, \cite{B}, \cite{J}, \cite{Ho}, \cite{BP}, \cite{D}, \cite{N},  \cite{DP}, \cite{DPM}, and \cite{BO}). \\
A generalized and interesting version of the classical CCP assumes (see, e.g., \cite{P}, \cite{FO}, \cite{FD}, \cite{Ad}, and \cite{Ro}), that the Dixie Cup company sells ice cream with a cardboard cover that has hidden on the underside a picture (``coupon'') of a sixties music band. In total there are $N$ different pictures and each one appears with probability $p_k$. Mr. and Mrs. Smith have one daughter and $(r-1)$ sons, all ice cream and sixties music addicts. The girl (she is the oldest) is the only one to buy ice cream. She tries to complete her collection. When she gets a new picture she puts it in her album, and when she gets a double, she gives it to her oldest brother, and when this one gets a double, he gives it to the remaining oldest brother, and so on. After having bought $T_N$ ice creams, the girl has completed her album while they remain $U_j^N$ unfilled places in the album of the $j$-th collector $j = 2, 3, \dots, r$ (that is the $(j-1)$-th brother). Obviously,
\begin{equation*}
1 \leq U_j^N \leq N
\qquad \text{and} \qquad
U_{j-1}^N \leq U_j^N,
\qquad
j = 2, 3, \dots ,
\end{equation*}
where, for completeness we have used the convention that $U_1^N = 0$. We shall refer to such a version of coupon collection as the \textbf{Generalized CCP (GCCP)}. In this paper we study the asymptotics of the expectation of the random variable $U_j^N$ as $N \to \infty$.
\subsection{The case of equal probabilities}
Naturally, the simplest case occurs when one
takes
\begin{equation*}
p_1 = \cdots = p_N = 1/N.
\end{equation*}
This case has been studied for quite a while. For $j = 2$, Pintacuda (see \cite{P}) used the martingale stopping theorem and proved that
\begin{equation}
E\left[U_2^N\right] = \sum_{m=1}^N \frac{1}{m} =: H_N
\label{00}
\end{equation}
($H_N$ is sometimes called the $N$-th harmonic number).\\\\
\textbf{Remark 1.} It is well known (see, e.g., \cite{F}) that
\begin{equation}
E[\,T_{N}\,]=N H_N.
\label{-1}
\end{equation}
Hence, for large $N$, when the main collector has completed her collection (notice that by (\ref{-1}) she will need in average $N\ln N + O(N)$ trials in order to succeed), the expected number of unfilled coupons in her oldest brother collection will be $\ln N + O(1)$.\\\\
Foata \textit{et} al., and Foata and Zeilberger (see \cite{FO} and \cite{FD}) using nonelementary mathematics, obtained recursive formulae for $E[U_j^N]$,
$j \geq 2$. Soonafter, Adler \textit{et} al. \cite{Ad} derived the same recursion, as well as a closed-form expression for $E[U_j^N]$ by using basic probability arguments (again for $j \geq 2$, while all $p_k$'s were considered equal). In particular, they proved that
\begin{equation}
E\left[U_2^N\right] = \sum_{m=1}^N \frac{1}{m},
\qquad
E\left[U_j^N\right]
= \sum_{m=1}^N \frac{E[U_{j-1}^m]}{m}
\quad \text{for }\, j \geq 3.
\label{r2}
\end{equation}
and also that
\begin{equation}
E\left[U_j^N\right] = \sum_{k=1}^N \binom{N}{k}\frac{(-1)^{k+1}}{k^{j-1}}
\qquad \text{for }\, j = 2, 3, \dots.
\label{r1}
\end{equation}
Foata \textit{et} al. (in \cite{FO}) called the quantites appearing in the recursion of (\ref{r2})
\textit{hyperharmonic numbers}. Notice that, for fixed $j$, detailed asymptotics for hyperharmonic numbers can be derived (e.g., via the assosiated generating functions). For example,
\begin{equation}
E\left[U_j^N\right] \sim \frac{(\ln N)^{j-1}}{(j-1)!}.\label{15}
\end{equation}
Furthermore, in the case of equal coupon probabilities, by exploiting the techniques of \cite{Ad} one can compute explicitly the variance $V[U_j^N]$ and its asymptotics as $N \to \infty$. In particular, for $j=2$ we get
\begin{equation}
V[U_2^N] = 4 \left(H_1 + \frac{H_2}{2} + \cdots + \frac{H_N}{N}\right) - 3 H_N
- H_N^2,
\label{Ross}
\end{equation}
(a slightly different, albeit equivalent, form of formula (\ref{Ross}) can be found in \cite{Ro1} and \cite{Ro})
hence
\begin{equation*}
V[U_2^N] \sim \ln^2 N
\qquad \text{as }\; N \to \infty.
\end{equation*}
%
%

\subsection{The case of unequal probabilities}
Let us now suppose that each coupon appears with probability $p_k$, with $\sum_{k=1}^{N}p_{k}=1,\,\,\,p_{k}>0$ for $k = 1,...,N$.
Then, Adler \textit{et} al. \cite{Ad}, (see also Ross \cite{Ro}) proved that,
when the main collector has a complete set, the expected number of unfilled coupons in each of her brothers' albums is obtained from
\begin{equation}
E[U_j^N]
= \sum_{k=1}^{N}\int_{0}^{\infty} p_{k}\,e^{-p_{k}t}\,\frac{\left(p_{k}t\right)^{j-1}}{\left(j-1\right)!}
\left[\prod_{i\neq k}\left(1-e^{-p_{i}t}\right)\right]dt,
\qquad
j \geq 2.
\label{r3}
\end{equation}
For example,
\begin{equation}
E\left[U_j^2\right] = 2 - p_1^j - p_2^j,
\qquad j \geq 2,
\label{EE1}
\end{equation}
and
\begin{equation}
E\left[U_2^3\right] = 3 + p_1^2 + p_2^2 + p_3^2
-\frac{p_1^2 + p_2^2}{(p_1 + p_2)^2} - \frac{p_1^2 + p_3^2}{(p_1 + p_3)^2}
-\frac{p_2^2 + p_3^2}{(p_2 + p_3)^2}.
\label{EE2}
\end{equation}
In addition, in \cite{Ad}, the authors derived two sets of general bounds for $E[\,U_{j}^{N}\,]$, as well as a simulation approach for estimating the
summands of (\ref{r3}).\\

\noindent \textbf{Conjecture.} For a fixed $N$, the average $E[U_j^N]$ of $U_j^N$ becomes maximum when all the $p_k$'s become equal (to $1/N$).\\
The results of the present paper support this conjecture.
\subsection{Large $N$ asymptotics}
When $N$ is large it is not clear at all what information one can obtain
from formula (\ref{r3}). For this reason there is a need to develop efficient
ways for deriving asymptotics for $E[\,U_{j}^{N}\,]$ as $N \rightarrow \infty $.\\
Let $\alpha =\{a_{k}\}_{k=1}^{\infty }$ be a sequence of strictly positive
numbers. Then, for each integer $N > 0$, one can create a probability measure
$\pi _N = \{p_1,...,p_N\}$ on the set of types $\{1,...,N\}$ by taking
\begin{equation}
p_k = \frac{a_k}{A_N},
\qquad \text{where}\quad
A_N = \sum_{k=1}^N a_k.
\label{B0}
\end{equation}
Notice that $p_k$ depends on $\alpha $ and $N$, thus, given $\alpha $, it
makes sense to consider the asymptotic behavior of $E[\,U_{j}^{N}\,] $ as $N\rightarrow \infty $. This approach for creating sequences of probability measures was first introduced in \cite{BP} and adopted in \cite{DP} and \cite{DPM}.\\
The sequence of measures $\pi_N$, $N = 2,3,...$, constucted from $\alpha$ via (\ref{B0}) has an interesting property:\\
For each $N = 2,3,...$ let
\begin{equation}
\Omega_N := \left\{\omega=\{\omega_j\}_{j=1}^{\infty}:\;
\omega_j = 1 \ \text{or } 2 \ \text{or } \cdots \text{or } N \right\},
\label{BB1}
\end{equation}
i.e. $\Omega_N = \{1, \dots, N \}^{\mathbb{N}}$, where
$\mathbb{N} = \{1, 2, \dots\}$. The sample space $\Omega_N$ describes the experiment of sampling $N$ coupons with replacement, indefinitely.

To each $m$-tuple $(k_1, \dots, k_m) \in \{1, \dots, N \}$ ($m \in \mathbb{N}$) we associate a cylinder subset of $\Omega_N$:
\begin{equation}
A_{k_1, \dots, k_m} := \{k_1\} \times \cdots \times \{k_m\}
\times \prod_{j = m+1}^{\infty} \{1, \dots, N \}.
\label{BB2}
\end{equation}
Then, $\pi_N = \{p_1,...,p_N\}$, as defined in (\ref{B0}), induces a set function (a probability) $P_N$ on these cylinder sets:
\begin{equation}
P_N\left\{ A_{k_1, \dots, k_m} \right\} := p_{k_1}^N \cdots p_{k_m}^N
\label{BB3}
\end{equation}
(here $N$ is a superscript indicating the dependence of $p_k$ on $N$). By  the Extension Theorem of Caratheodory, $P_N$ extends to a complete probability measure (which we also denote by $P_N$) on $(\Omega_N, \cal{F}_N)$, where $\cal{F}_N$ is the completion of the $\sigma$-algebra generated by the cylinder sets.\\
If $L_{N+1}$ is the subset of $\Omega_{N+1}$ defined by
\begin{equation}
L_{N+1} := \left\{\omega=\{\omega_j\}_{j=1}^{\infty}:\;
\omega_j = N+1 \ \text{for all but finitely many } j \right\},
\label{BB4}
\end{equation}
then it is clear that
\begin{equation}
P_{N+1}\left\{ L_{N+1} \right\} = 0.
\label{BB5}
\end{equation}
For a sequence
$\omega = \{\omega_j\}_{j=1}^{\infty} \in \Omega_{N+1} \setminus L_{N+1}$
we denote by $\iota_{N+1}(\omega)$ the sequence (in $\Omega_N$) which is obtained
from $\omega$ by deleting all terms $\omega_j$ such that $\omega_j = N+1$
(notice that $\iota_{N+1}(\omega)$ is a sequence, i.e. it has infinitely many terms, since $\omega \notin L_{N+1}$).\\
To be more precise, given $\omega = \{\omega_j\}_{j=1}^{\infty}$ let us consider the set of indices
\begin{equation*}
\{j_n \in \mathbb{N} :\; \omega_{j_n} \ne N+1, \quad j_n < j_{n+1}\}.
\end{equation*}
Then, $\iota_{N+1}(\omega) = \{\omega_{j_n}\}_{n=1}^{\infty}$. We call $\iota_{N+1}(\omega)$ the $(N+1)$-\textit{contraction} of
$\omega \in \Omega_{N+1} \setminus L_{N+1}$. For example, if
\begin{equation*}
\omega = (3, 1, 2, 3, 3, 2, 1, 2, 2, 3, 1, 1, \dots) \in \Omega_3 \setminus L_3,
\end{equation*}
then
\begin{equation*}
\iota_3(\omega) = (1, 2, 2, 1, 2, 2, 1, 1, \dots) \in \Omega_2.
\end{equation*}
Since $\Omega_N \subset (\Omega_{N+1} \setminus L_{N+1})$, we have that
$\iota_{N+1} : (\Omega_{N+1} \setminus L_{N+1}) \to \Omega_N$ is a
``projection'' and, via $\iota_{N+1}$, $P_{N+1}$ induces a measure $\tilde{P}_N$ on
$\cal{F}_N$ given by
\begin{equation}
\tilde{P}_N\{A\} := P_{N+1}\{\iota_{N+1}^{-1}(A)\}.
\label{BB5b}
\end{equation}

\noindent \textbf{Proposition 1.} The measures $P_N$ and $\tilde{P}_N$ coincide, namely
\begin{equation*}
\tilde{P}_N\{A\} = P_N\{A\}
\qquad \text{for all }\; A \in \cal{F}_N.
\end{equation*}
The proof is straightforward. \\\\
\textbf{Remark 2.} For $\omega = \{\omega_j\}_{j=1}^{\infty} \in \Omega_{N+1}
\setminus L_{N+1}$ we have that $T_{N+1}(\omega) =$ the smallest $k$ such that $(\omega_1, \dots, \omega_k)$ contains all elements of $\{1, \dots, N+1\}$. If we set $\tilde{T}_N(\omega) := T_N(\iota_{N+1}(\omega))$, then, obviously $\tilde{T}_N(\omega) \leq T_{N+1}(\omega)$. Also, it is easy to see that
$\tilde{T}_N$ and $T_N$ have the same distribution. Therefore the sequence $\{T_N\}_{N = 2}^{\infty}$ is \textit{stochastically increasing}, i.e. $P\{T_{N+1} \geq k\} \geq P\{T_N \geq k\}$ for all $k$ and all $N \geq 2$.\\
This is not true however, for the sequence $\{U^N_j\}_{N = 2}^{\infty}$. For instance, let $\alpha = (1, 1, \epsilon, \dots)$ where $\epsilon > 0$. Then, for $N = 2$ formula (\ref{B0}) gives that $p_1 = p_2 = 1 / 2$ and hence, by (\ref{EE1}) (or (\ref{00})) we get that $E\left[U_2^2\right] = 3 / 2$. On the other hand, if $N = 3$, then (\ref{B0}) gives $p_1 = p_2 = 1 / (2 + \epsilon)$ and $p_3 = \epsilon / (2 + \epsilon)$. It, then, follows from (\ref{EE2}) that
$E\left[U_3^2\right]$ can be made as close to $1$ as we wish, by taking $\epsilon$ sufficiently close to $0$ (as it is intuitively expected, since one coupon is very ``rare'' compare to the others). In particular, we can have
$E\left[U_3^2\right] < E\left[U_2^2\right]$, for a given sequence $\alpha$.\\
If $\alpha$ is the constant sequence whose terms are (all) equal to $1$ (i.e. in the case of equal coupon probabilities), then (\ref{r2}) implies that $E[U_j^N]$ increases with $N$. We conjecture that in this case we also have that, for any $j \geq 2$, $U_j^N$ is stochastically increasing with $N$. \\\\
After the above deviation we return to our main issue, namely the asymptotics of $E[U_j^N]$. Inspired by (\ref{r3}) we introduce the notation
\begin{equation}
I_N(\alpha; j) := \sum_{k=1}^{N} \int_0^{\infty} a_k \,
e^{-a_k t} \, \frac{(a_k t)^{j-1}}{(j-1)!}
\prod_{i \neq k} \left(1 - e^{-a_i t}\right)\,dt,
\qquad j \geq 2.
\label{4}
\end{equation}
If $s \alpha := \{sa_k\}_{k=1}^{\infty}$, (\ref{4}) gives immediately that
\begin{equation}
I_{N}(s\alpha ) =  I_N(\alpha)
\label{4a}
\end{equation}
and hence, in view of (\ref{r3}) and (\ref{B0})
\begin{equation}
E\left[U_j^N\right] = I_N(\alpha; j).
\label{5}
\end{equation}
The rest of the paper is organized as follows. In Section 2 we consider classes of decaying sequences
$\alpha$ such that $a_k \rightarrow 0$. Here the computations are quite involved. We present the main result in Theorem 1 (of Subsection 2.3). In particular, the first three terms of the asymptotic expansion of $E[\,U_{j}^{N}\,]$ are determined (as $N\rightarrow \infty$). It is notable that the \textit{generalized Zipf law} falls in this category. The method of proving Theorem 1 is based on a ``brute force'' technique reminding the technique initiated in \cite{BP} and exploited in \cite{DP}.
In Section 3 we derive the leading behavior of $E[\,U_{j}^{N}\,]$ for a large class of sequences $\alpha$, such that $a_k \rightarrow \infty$. In this case we often get that $E[U_j^N]$ approaches a finite limit as $N \rightarrow \infty$. Various examples are exhibited. In particular, we cover some important families of coupon probabilities (e.g. polynomial and exponential).

\section{Decaying sequences}
Inspired by \cite{DP} we consider sequences $\alpha =\{a_{k}\}_{k=1}^{\infty }$ of the form
\begin{equation}
a_{k}=\frac{1}{f(k)}, \label{aj}
\end{equation}
where
\begin{equation}
f(x) > 0
\quad \text{and}\quad
f'(x) > 0,
\qquad x > 0,
\label{C1a}
\end{equation}
and furthermore we assume that $f(x)$ possesses three derivatives and satisfies
the following conditions as $x \rightarrow \infty$:
\begin{equation*}
\;\;\;\;\;\;\;\;\;\;\;\;\text{(i) }
f(x)\rightarrow \infty,
\;\;\;\;\;\;\;\;\;\;\;\;\;\;\;\;\;\;\;\;\;\;\;\;\;\;\; \text{(ii) }
\frac{f^{\prime }(x)}{f(x)}=o\left(\frac{1}{\ln^{p}x}\right)\,\,\text{for any } p > 0,
\end{equation*}
\begin{equation}
\text{(iii) }
\frac{f^{\prime \prime}(x)/f^{\prime }(x)}{f'(x)/f(x)} = O\left(1\right),
\;\;\;\;\;\;\;\;\;\;\;\text{(iv) }
\frac{f^{\prime \prime\prime}(x)\;f(x)^{2}}{ f^{\prime }(x)^{3}} = O\left(1\right)
\label{C1}
\end{equation}
(in \cite{DP} the conditions on $f(x)$ were slightly weaker). These conditions are satisfied by a variety of commonly
used functions. For example,
\begin{equation*}
f(x) = x^p (\ln x)^q, \quad p > 0,\ q \in \mathbb{R},\qquad \qquad
f(x) = \exp(x^{r}),\quad 0 < r < 1,
\end{equation*}
or various convex  combinations of products of such functions.\\\\
\textbf{Remark 3.} Condition (ii) of (\ref{C1}) implies $f^{\prime }(x)/f(x)\rightarrow 0$. Thus,
\begin{equation}
\lim_{x \rightarrow \infty}\frac{f(x+1)}{f(x)} = 1
\label{LIMIT}
\end{equation}
(this can be justified, e.g., by applying the Mean Value Theorem to the function $g(x) = \ln f(x)$
on the interval $[x, x+1]$).\\\\
For typographical convenience we set
\begin{equation}
F(x) := f(x) \ln\left(\frac{f(x)}{f'(x)}\right)
\label{F}
\end{equation}
(notice that (\ref{C1a}) and (ii) of (\ref{C1}) imply that $F(x) > 0$ for $x$ sufficiently large). Starting from (\ref{4}), we substitute $t = F(N) s$ in the
integral and rewrite $I_{N}(\alpha;j )$ as\\\\
\begin{equation}
I_{N}(\alpha;j)=\frac{1}{\left(j-1\right)!}\,F(N)^{j}\left[I_{N}^{1} (\alpha;j )+I_{N}^{2} (\alpha;j )\right],\label{20}
\end{equation}
where
\begin{equation}
I_{N}^{1} (\alpha;j ): = \sum_{k=1}^{N}f(k)^{-j}\int_0^1 e^{-\frac{F(N)}{f(k)}s}\exp \left[ \sum_{i\neq k}\ln \left( 1-e^{-\frac{F(N)}{f(i)}s}\right) \right]\,s^{j-1} ds
\label{I1}
\end{equation}
and
\begin{equation}
I_{N}^{2} (\alpha;j ):= \sum_{k=1}^{N}f(k)^{-j}\int_1^{\infty} e^{-\frac{F(N)}{f(k)}s}\exp \left[ \sum_{i\neq k}\ln \left( 1-e^{-\frac{F(N)}{f(i)}s}\right) \right]\,s^{j-1} ds.
\label{I2}
\end{equation}
In order to analyze deeper the above quantities we need the following lemma.\\\\
\textbf{Lemma 1.} Set
\begin{equation}
J_{m}(N):=\int_{1}^{N} f(x)^{m}e^{-\frac{F(N)}{f(x)}s}dx ,\quad m\geq0.
\label{I}
\end{equation}
Then, under (\ref{C1}) and (\ref{F}), we have, as $N\rightarrow \infty$,
\begin{align}
J_{m}(N)&=\frac{f(N)^{m+2}}{sF(N)f^{\prime}(N)}e^{-\frac{F(N)}{f(N)}s} \nonumber\\
&+\omega(N)\;\frac{f(N)^{m+3}}{s^{2}F(N)^{2}f^{\prime}(N)}e^{-\frac{F(N)}{f(N)}s}\left[1+O\left(\frac{f(N)}{F(N)}\right)\right], \label{IIa}
\end{align}
uniformly in $s\in [s_{0},\infty)$, for any $s_{0} > 0 $,
where
\begin{equation}
\omega(N) := -2+\frac{f^{\prime \prime
}(N)/f^{\prime }(N)}{f^{\prime }(N)/f(N)}.
\label{a}
\end{equation}

For the proof see \cite{DP}. Notice that the condition (iii) of (\ref{C1}) says that $\omega(N) = O(1)$ as $N \to \infty$. \\\\
\textbf{Remark 4.} It is straightforward to check that Lemma 1 is still valid when $m$ is a negative integer.
\subsection{The integral $I_{N}^{1} (\alpha;j )$}
Regarding the quantity of (\ref{I1}), given $\varepsilon \in (0,1/2)$ we have
\begin{equation}
I_{N}^{1} (\alpha;j )=I_{N}^{11} (\alpha;j )+I_{N}^{12} (\alpha;j ),\label{-11}
\end{equation}
where
\begin{equation}
I_{N}^{11} (\alpha;j ): = \sum_{k=1}^{N}f(k)^{-j}\int_0^{1 - \varepsilon} e^{-\frac{F(N)}{f(k)}s}\exp \left[ \sum_{i\neq k}\ln \left( 1-e^{-\frac{F(N)}{f(i)}s}\right) \right]\,s^{j-1} ds \label{I1a}
\end{equation}
and
\begin{equation}
I_{N}^{12} (\alpha;j ):=\sum_{k=1}^{N}f(k)^{-j}\int_{1 - \varepsilon}^1 e^{-\frac{F(N)}{f(k)}s}\exp \left[ \sum_{i\neq k}\ln \left( 1-e^{-\frac{F(N)}{f(i)}s}\right) \right]\,s^{j-1} ds.
\label{I1b}
\end{equation}
The quantity of (\ref{I1a}) becomes (in view of (\ref{C1a}))
\begin{align}
I_{N}^{11} (\alpha;j )=&\sum_{k=1}^{N}f(k)^{-j}\int_0^{1 - \varepsilon} \frac{e^{-\frac{F(N)}{f(k)}s}}{1-e^{-\frac{F(N)}{f(k)}s}}\exp \left[ \sum_{i=1}^{N}\ln \left( 1-e^{-\frac{F(N)}{f(i)}s}\right) \right]\,s^{j-1} ds \nonumber\\
<&\sum_{k=1}^{N}f(k)^{-j}\int_0^{1 - \varepsilon} \frac{e^{-\frac{F(N)}{f(k)}s}}{1-e^{-\frac{F(N)}{f(k)}s}}
\exp \left[ \sum_{i=1}^{N}\ln \left( 1-e^{-\frac{F(N)}{f(i)}\left(1-\varepsilon\right)}\right) \right]\,s^{j-1} ds\nonumber\\
<& \exp \left(-\sum_{i=1}^N e^{-\frac{F(N)}{f(i)} \left(1 - \varepsilon \right)} \right)
\sum_{k=1}^{N}f(k)^{-j}\int_0^{1 - \varepsilon} \frac{e^{-\frac{F(N)}{f(k)}s}}{1-e^{-\frac{F(N)}{f(k)}s}}
\,s^{j-1} ds,
\label{21}
\end{align}
since $\ln(1-x) < -x$, for $0 < x < 1$. Now, $f$ is increasing, hence, from the comparison of sums and integrals we have
\begin{align}
\int_{1}^{N}e^{-\frac{F(N)}{f(x)}s}dx&\leq \sum_{j=1}^{N}e^{-\frac{F(N)}{f(j)}s} \nonumber\\
&\leq \int_{1}^{N+1}e^{-\frac{F(N)}{f(x)}s}dx \nonumber\\
&\leq \int_{1}^{N}e^{-\frac{F(N)
}{f(x)}s}dx+e^{-\frac{F(N)}{f(N+1)}s}. \label{0}
\end{align}
Using the above comparison, (\ref{LIMIT}), and applying Lemma 1, for $m = 0$, one arrives at
\begin{align}
&\exp \left(-\sum_{i=1}^N e^{-\frac{F(N)}{f(i)} \left(1 - \varepsilon \right)} \right)\leq\exp \left(-\int_{1}^{N} e^{-\frac{F(N)}{f(x)}\left(1-\varepsilon \right)}dx \right) \nonumber\\
&<\exp\left\{-\frac{f(N)^2}{(1 - \varepsilon) F(N) f'(N)} e^{-\frac{F(N)}{f(N)}(1 - \varepsilon)}\left(1 + M_1\;\frac{f(N)}{F(N)} \right) \right\}\nonumber\\
&=\exp\left\{-\frac{1}{1-\varepsilon} \cdot \frac{\left(\frac{f(N)}{f'(N)}\right)^{\varepsilon}}{\ln\left( \frac{f(N)}{f'(N)} \right)}
\left(1 + \frac{M_1}{\ln\left(\frac{f(N)}{f^{\prime}(N)}\right)}\right)\right\},\label{b2}
\end{align}
where $M_1$ is a positive constant and we have used (\ref{F}), i.e. the definition of $F$. On the other hand, if we set
\begin{equation*}
J_{N}(k;\varepsilon):=\int_0^{1 - \varepsilon} \frac{e^{-\frac{F(N)}{f(k)}s}}{1-e^{-\frac{F(N)}{f(k)}s}}
\,s^{j-1} ds,
\end{equation*}
then the scaling $u = F(N) f(k)^{-1}s$, via the definition of $F$ and the monotonicity (and positivity) of $f$, yields
\begin{align*}
J_{N}(k;\varepsilon)=&\left[\frac{f(k)}{F(N)}\right]^{j}\int_0^{\left(1 - \varepsilon\right)\frac{F(N)}{f(k)}}\frac{e^{-u}}{1-e^{-u}}\;u^{j-1}du \nonumber\\
<&\left[\frac{f(k)}{F(N)}\right]^{j}\int_0^{\infty}\frac{e^{-u}}{1-e^{-u}}\;u^{j-1}du.
\end{align*}
Since $j\geq 2$, the last integral is equal to a positive constant. In particular, for $j=2,$ it converges to $\pi^{2}/6.$ In general, it is
not difficult to check that it converges to
\begin{equation*}
M_j := \Gamma(j)\,\zeta(j),
\qquad
j = 2, 3,\dots,
\end{equation*}
where $\Gamma(\cdot),\,\zeta(\cdot)$ denote the gamma function and the Riemann zeta function respectively. Hence,
\begin{equation}
J_{N}(k;\varepsilon)\leq M_j\, \left[\frac{f(k)}{F(N)}\right]^{j}.\label{22}
\end{equation}
In view of (\ref{b2}) and (\ref{22}), (\ref{21}) yields
\begin{align}
I_{N}^{11}(\alpha;j)\leq M_{j}\frac{N}{F(N)^{j}}\,\exp&\left\{-\frac{1}{1-\varepsilon} \cdot \frac{\left(\frac{f(N)}{f'(N)}\right)^{\varepsilon}}{\ln\left( \frac{f(N)}{f'(N)} \right)}\left(1 + M_1\;\frac{1}{\ln\left(\frac{f(N)}{f^{\prime}(N)}\right)}\right)\right\}.
\label{123}
\end{align}
Since $\varepsilon \in (0,1/2)$, we claim that (\ref{123}) implies
\begin{equation}
I_{N}^{11}(\alpha;j)
<< \frac{\ln^{2}\left(\ln\left(\frac{f(N)}{f^{\prime}(N)}\right)\right)}{f(N)^{j}\ln^{3}\left(\frac{f(N)}{f^{\prime}(N)}\right)},
\qquad j \geq 2,
\label{b3}
\end{equation}
where $a_N << b_N$ means that $a_N / b_N \to 0$ as $N \to \infty$.
To check the validity of (\ref{b3}) one observes that it suffices to show that
\begin{equation}
N<< \exp\left[\frac{1}{1-\varepsilon} \cdot \frac{\left(\frac{f(N)}{f'(N)}\right)^{\varepsilon}}{\ln\left( \frac{f(N)}{f'(N)} \right)}\right]
\times\left[\frac{\ln^{2}\left(\ln\left(\frac{f(N)}{f^{\prime}(N)}\right)\right)}{\ln^{3-j}\left(\frac{f(N)}{f^{\prime}(N)}\right)}\right],\label{01}
\end{equation}
which follows by taking logarithms and using condition (ii) of (\ref{C1}).\\
Our next task is to compute a few terms of the asymptotic expansion of the term $I_{N}^{12} (\alpha;j )$ defined in (\ref{I1b}). For convenience we set
\begin{equation}
B_{k}(N;s) := \sum_{i=1}^N\ln \left( 1 - e^{-\frac{F(N)}{f(i)} s} \right)-\ln \left( 1 - e^{-\frac{F(N)}{f(k)} s} \right).
\label{AsN}
\end{equation}
Since
\begin{equation*}
\frac{F(N)}{f(i)}\rightarrow \infty \qquad \text{as }\ N\rightarrow
\infty,
\end{equation*}
and $\ln (1-x) = -x + O(x^2)$ as $x\rightarrow 0$, we have (as long as $s\geq s_{0} > 0$)
\begin{equation}
B_{k}(N;s)=\sum_{i=1}^{N} \left[-e^{-\frac{F(N)}{f(i)}s}+O\left(e^{-\frac{2F(N)}{f(i)}s}\right)\right]+e^{-\frac{F(N)}{f(k)} s}+O\left(e^{-\frac{2F(N)}{f(k)} s}\right).\label{O1}
\end{equation}
From the comparison of sums and integrals, i.e. (\ref{0}), (\ref{O1}) yields
\begin{align*}
B_{k}(N;s)=-&\left[\int_{1}^{N}e^{-\frac{F(N)}{f(x)}s}dx+O\left(e^{-\frac{F(N)}{f(N+1)}s}\right)\right] \nonumber\\
+&\sum_{i=1}^{N}O\left(e^{-\frac{2F(N)}{f(i)}s}\right)+e^{-\frac{F(N)}{f(k)} s}+O\left(e^{-\frac{2F(N)}{f(k)} s}\right).
\end{align*}
The above formula together with Lemma 1, for $m = 0$, give
\begin{equation*}
B_{k}(N;s)=-\frac{f(N)^{2}}{sF(N)f^{\prime}(N)}e^{-\frac{F(N)}{f(N)}s}-\omega(N)\;\frac{f(N)^{3}}{s^{2}F(N)^{2}f^{\prime}(N)}e^{-\frac{F(N)}{f(N)}s}\left[1+O\left(\frac{f(N)}{F(N)}\right)\right]
\end{equation*}
\begin{equation*}
\;\;\;\;\;\;\;\;\;\;+O\left(e^{-\frac{F(N)}{f(N+1)}s}+Ne^{-\frac{2F(N)}{f(N)}s}\right)+e^{-\frac{F(N)}{f(k)} s}+O\left(e^{-\frac{2F(N)}{f(k)} s}\right).
\end{equation*}
Using (\ref{LIMIT}) the above yields
\begin{equation*}
B_{k}(N;s) = -\frac{f(N)^{2}}{sF(N)f'(N)}e^{-\frac{F(N)}{f(N)} s}
-\omega(N)\;\frac{f(N)^{3}}{s^{2}F(N)^{2}f'(N)}e^{-\frac{F(N)}{f(N)} s}\left[1+O\left(\frac{f(N)}{F(N)}\right)\right],
\end{equation*}
independent of $k$. Hence,
\begin{align*}
I_{N}^{12} (\alpha;j )=&\sum_{k=1}^{N}f(k)^{-j} \int_{1-\varepsilon}^{1} e^{B(N;s)}\, e^{-\frac{F(N)}{f(k)}s}\,s^{j-1} ds\nonumber\\
=&\sum_{k=1}^{N}f(k)^{-j}\int_{1-\varepsilon}^1 e^{-\frac{F(N)}{f(k)}s}\,s^{j-1}\times\exp\left[-\frac{f(N)^{2}}{sF(N)f'(N)}e^{-\frac{F(N)}{f(N)}s}\right.\nonumber\\
&\left.\;\;\;\;\;\;-\omega(N)\;\frac{f(N)^3}{s^2 F(N)^{2}f'(N)}e^{-\frac{F(N)}{f(N)}s}\left[1 + O\left(\frac{f(N)}{F(N)}\right)\right]\right]ds
\end{align*}
as $N \rightarrow \infty$. Using the definition of $F$, namely (\ref{F}), and substituting $s = 1 - t$, the above expression becomes
\begin{align*}
I_N^{12} (\alpha;j) = \sum_{k=1}^{N}f(k)^{-j}&\int_0^{\varepsilon}\left(1-t\right)^{j-1}\exp\left[-\left(1-t\right)\frac{f(N)}{f(k)}\ln\left(\frac{f(N)}{f'(N)}\right)\right] \nonumber\\
&\;\;\;\times
\exp\left\{-\frac{1}{1-t}\;\frac{\left(\frac{f(N)}{f'(N)}\right)^{t}}{\ln\left(\frac{f(N)}{f'(N)}\right)}
-\omega(N)\frac{1}{\left(1-t\right)^{2}}\;\frac{\left(\frac{f(N)}{f^{\prime}(N)}\right)^{t}}{\ln\left(\frac{f(N)}{f^{\prime}(N)}\right)^2}\right.\nonumber\\
&\left.\;\;\;\;\;\;\;\;\;\;\;\;\;\;\;\;\;\;\;\;\;\;\;\;\;\;\;\;\;\;\;\;\;\;\;\;\;\;\;\;\;\;\;\;\;\times\left[1 + O\left(\frac{1}{\ln\left(\frac{f(N)}{f'(N)}\right)}\right)\right]\right\} dt.
\end{align*}
For typographical convenience we set
\begin{equation}
A := \frac{f(N)}{f'(N)}
\label{A}
\end{equation}
(notice that $A \rightarrow \infty$ as $N \rightarrow \infty$). Then $I_N^{12}$ can be expressed as
\begin{align}
I_{N}^{12} (\alpha;j)=\int_{0}^{\varepsilon}&\left[\sum_{k=1}^{N}f(k)^{-j}\,\exp\left(-\frac{1}{f(k)}f(N)\left(1-t\right)\ln A\right)\right] \nonumber\\
&\times \left(1-t\right)^{j-1}
\exp\left\{-\;\frac{A^{t}}{\ln A}\left(\sum_{n=0}^{\infty}t^{n}\right)-\omega(N)\frac{A^{t}}{\ln^{2}A}\;\left(\sum_{n=1}^{\infty}nt^{n-1}\right)\right.\nonumber\\
&\left.\;\;\;\;\;\;\;\;\;\;\;\;\;\;\;\;\;\;\;\;\;\;\;\;\;\;\;\;\;\;\;\;\;\;\;\;\;\;\;\;\;\;\;\;\;\;\;\;\;\;\;\;\;\;\;\times\left[1+O\left(\frac{1}{\ln A}\right)\right]\right\}dt.
\label{eint1}
\end{align}
Set
\begin{equation*}
S_j(N;t):= \sum_{k=1}^{N}f(k)^{-j}\,e^{-\frac{f(N)}{f(k)} \left(1-t\right)\ln A},\,\,\,\,\,\,\,\,\,\,\,0\leq t \leq \varepsilon<1/2.
\end{equation*}
and
\begin{equation*}
g_{j}(x):=f(x)^{-j}\,\exp\left(-\frac{f(N)}{f(x)} \left(1-t\right)\ln A\right),\;\;0\leq x\leq N, \;\;0\leq t \leq \varepsilon<1/2.
\end{equation*}
It is easy to check that under conditions (\ref{C1}), $g$ is increasing for sufficiently large  $N$. Thus, it follows from the
comparison of sums and integrals that
\begin{equation}
S_j(N;t)=K_j(N;t)+O\left(\frac{e^{-\frac{f(N)}{f(N+1)}\left(1-t\right)\ln A}}{f(N+1)^{j}}\right),\label{S}
\end{equation}
where
\begin{equation}
K_{j}(N;t):=\int_{1}^{N}f(x)^{-j}\,e^{-\frac{f(N)}{f(x)}\left(1-t\right)\ln A}dx.\label{K}
\end{equation}
By using Lemma 1 and Remarks 4 and 5 (as long as $1-t>1/2$), we get as $N\rightarrow \infty$,
\begin{align}
S_{j}(N;t)&=\frac{A^{t-1}}{\left(1-t\right)\ln A}\cdot\frac{1}{f(N)^{j-1}f^{\prime}(N)} \nonumber\\
&+\omega(N)\;\frac{A^{t-1}}{\left(1-t\right)^{2}\ln^{2} A}\cdot\frac{1}{f(N)^{j-1}f^{\prime}(N)}\left[1+O\left(\frac{1}{\ln A}\right)\right] \nonumber\\
&+O\left(\frac{A^{t-1}}{f(N+1)^{j}}\right).\label{23}
\end{align}
In view of (\ref{23}) and (\ref{A}), (\ref{eint1}) yields
\begin{align*}
I_{N}^{12} (\alpha;j )&=\frac{1}{ f(N)^{j-1}f^{\prime}(N)\,A\,\ln A} \nonumber\\
&\times\int_{0}^{\varepsilon}\left[\frac{A^{t}}{1-t}+\omega(N)\;\frac{A^{t}}{\left(1-t\right)^{2}\ln A}\left[1+O\left(\frac{1}{\ln A}\right)\right]+O\left( A^{t-1}\ln A\right)\right] \nonumber\\
&\;\;\;\;\;\;\;\;\;\;\times \left(1-t\right)^{j-1}\exp\left[-\;\frac{A^{t}}{\ln A}\left(\sum_{n=0}^{\infty}t^{n}\right)-\omega(N)\frac{A^{t}}{\ln^{2}A}\;\left(\sum_{n=1}^{\infty}nt^{n-1}\right)\right.\nonumber\\
&\left.\;\;\;\;\;\;\;\;\;\;\;\;\;\;\;\;\;\;\;\;\;\;\;\;\;\;\;\;\;\;\;\;\;\;\;\;\;\;\;\;\;\;\;\;\;\;\;\;\;\;\;\;\;\;\;\;\;\;\;\;\;\;\;\;\;\;\;\times\left[1+O\left(\frac{1}{\ln A}\right)\right]\right]dt.
\end{align*}
Substituting $u = A^t / \ln A$ in the integral above, we get (in view of (\ref{A}))
\begin{align*}
I_{N}^{12} (\alpha;j )&
= \frac{1}{f(N)^{j}\ln A}\int_{1/\ln A}^{A^{\varepsilon}/\ln A}
\left\{\left[\frac{1}{1-\frac{\ln u}{\ln A}-\frac{\ln(\ln A)}{\ln A}}\right.\right.\nonumber\\
&\left.\left.\qquad\qquad\quad+\frac{\omega(N)}{\ln A}\;\frac{1}{\left(1-\frac{\ln u}{\ln A}-\frac{\ln(\ln A)}{\ln A}\right)^{2}}\left[1+O\left(\frac{1}{\ln A}\right)\right]+O\left(\frac{\ln A}{A}\right)\right]
\right.
\\
&
\left.
\qquad\qquad\quad\times \left(1-\frac{\ln u}{\ln A}-\frac{\ln(\ln A)}{\ln A}\right)^{j-1}\right.\nonumber\\
&\left.\times\exp\left[-\frac{u}{1-\frac{\ln u}{\ln A}-\frac{\ln(\ln A)}{\ln A}}-\frac{\omega(N)}{\ln A}\;\frac{u}{\left(1-\frac{\ln u}{\ln A}-\frac{\ln(\ln A)}{\ln A}\right)^{2}}\left[1+O\left(\frac{1}{\ln A}\right)\right]\right]\,\right\}du.
\end{align*}
If we set
\begin{equation}
 \delta := \frac{1}{\ln A}=\frac{1}{\ln \left(\frac{f(N)}{f^{\prime}(N)}\right)}=\frac{f(N)}{F(N)}\label{delta}
\end{equation}
(hence, $A\rightarrow \infty$ implies $\delta \rightarrow 0^+$), the above integral becomes
\begin{align}
I_{N}^{12} (\alpha;j )&=\frac{\delta}{ f(N)^{j}}\int_{\delta}^{\delta \exp\left(\varepsilon/\delta \right)}\left\{\left[\left(1-\delta \ln u+\delta \ln \delta\right)^{j-2}\right.\right.\nonumber\\
&\left.\left.\qquad\qquad+\omega(N)\delta\left(1-\delta \ln u+\delta \ln \delta\right)^{j-3}\left(1+O\left(\delta\right)\right)+O\left(\frac{e^{-1/ \delta}}{\delta}\right)\right]\right.\nonumber\\
&\left.\qquad\qquad\times \exp\left[-\frac{u}{1-\delta \ln u+\delta \ln \delta}-\frac{\omega(N)u \delta}{\left(1-\delta \ln u+\delta \ln \delta\right)^{2}}\left(1+O\left(\delta\right)\right)\right]\,\right\}du\label{eint5}
\end{align}
We split the integral of (\ref{eint5}) as:
\begin{equation}
\int_{\delta}^{\delta \exp\left(\varepsilon/\delta \right)}=\int_{\delta}^{1/\sqrt{\delta}}+\int_{1/\sqrt{\delta}}^{\delta \exp(\varepsilon/\delta)}.\label{aw}
\end{equation}
The second integral of (\ref{aw}) can be bounded as follows:
\begin{align*}
&\int_{1/\sqrt{\delta}}^{\delta \exp(\varepsilon/\delta)}\left\{\left[\left(1-\delta \ln u+\delta \ln \delta\right)^{j-2}\right.\right.\nonumber\\
&\left.\left.\qquad\qquad+\omega(N)\delta\left(1-\delta \ln u+\delta \ln \delta\right)^{j-3}\left(1+O\left(\delta\right)\right)+O\left(\frac{e^{-1/ \delta}}{\delta}\right)\right]\right.\nonumber\\
&\left.\qquad\qquad\times \exp\left[-\frac{u}{1-\delta \ln u+\delta \ln \delta}-\frac{\omega(N)u \delta}{\left(1-\delta \ln u+\delta \ln \delta\right)^{2}}\left(1+O\left(\delta\right)\right)\right]\,\right\}du
\end{align*}
\begin{align}
=&\int_{1/\sqrt{\delta}}^{\delta \exp(\varepsilon/\delta)}\left\{\left(1-\delta \ln u + \delta \ln \delta\right)^{j-2}\left[1+\frac{\omega(N)}{1-\delta \ln u+\delta \ln \delta}\delta \left(1+O\left(\delta\right)\right)+O\left(\frac{e^{-1/ \delta}}{\delta}\right)\right]\right.\nonumber\\
&\left.\;\;\;\;\;\;\;\;\;\;\;\;\;\;\;\;\;\;\;\;\times \exp\left[-\frac{u}{1-\delta \ln u + \delta \ln \delta}\left[ 1 + \frac{\omega(N) }{1-\delta \ln u + \delta \ln \delta}\;\delta\left(1+O\left(\delta\right)\right)\right]\right]\right\}\,du\nonumber\\
&\leq M \int_{1/\sqrt{\delta}}^{\infty} e^{-u}\,du=O\left( e^{-1/\sqrt{\delta}}\right),\label{eint6}
\end{align}
for some positive constant $M$ (since, $0<\varepsilon<1/2$). The first integral of (\ref{aw}) is
\begin{align*}
K_1(\delta) := &\int_{\delta}^{1/\sqrt{\delta}}\left\{\left[\left(1-\delta \ln u+\delta \ln \delta\right)^{j-2}\right.\right.\nonumber\\
&\left.\left.\qquad\qquad+\omega(N)\delta\left(1-\delta \ln u+\delta \ln \delta\right)^{j-3}\left(1+O\left(\delta\right)\right)+O\left(\frac{e^{-1/ \delta}}{\delta}\right)\right]\right.\nonumber\\
&\left.\qquad\qquad\times \exp\left(-\frac{u}{1-\delta \ln u+\delta \ln \delta}-\frac{\omega(N)u \delta}{\left(1-\delta \ln u+\delta \ln \delta\right)^{2}}\left[1 + O\left(\delta\right)\right]\right)\,\right\}du
\end{align*}
Since $\left( 1 - x \right)^{-2} =\sum_{n=1}^{\infty} n x^{n-1}$ for $\left|x\right| < 1$, we have
\begin{equation*}
K_1(\delta) = \int_{\delta}^{1/\sqrt{\delta}}
\left[\left(1-\delta \ln\frac{u}{\delta}\right)^{j-2}+\omega(N)\delta\left(1-\delta \ln\frac{u}{\delta}\right)^{j-3}\left(1+O\left(\delta\right)\right)+O\left(\frac{e^{-1/ \delta}}{\delta}\right)\right]
\end{equation*}
\begin{equation*}
\;\;\;\;\;\;\;\;\;\;\;\;\;\;\;\;\;\;\;\;\;\;\;\;\;\;
\times\exp\left[-u\sum_{n=0}^{\infty} \left( \delta \ln\frac{u}{\delta} \right)^n - u \; \omega(N)\delta \left(1 + O\left(\delta\right)\right) \sum_{n=1}^{\infty} n \left( \delta \ln\frac{u}{\delta}\right)^{n-1}\right]\,du.
\end{equation*}
We use the binomial theorem to expand the quantities $\left(1-\delta \ln\frac{u}{\delta}\right)^{j-2}$ and $\left(1-\delta \ln\frac{u}{\delta}\right)^{j-3}$ and get
\begin{align*}
K_1(\delta)=&\int_{\delta}^{1/\sqrt{\delta}}\left\{\left[\left(1-\left(j-2\right)\delta \ln\frac{u}{\delta}+O\left(\delta^{2} \ln^{2}\frac{u}{\delta}\right)\right)\right.\right.\nonumber\\
&\left.\left.\qquad\quad+\omega(N)\,\delta\,\left(1+O\left(\delta\right)\right)\left(1+O\left(\delta \ln\frac{u}{\delta}\right)\right)+O\left(\frac{e^{-1/ \delta}}{\delta}\right)\right]\right.\nonumber\\
&\left.\times\,e^{-u}\;\exp\left(-u\sum_{n=1}^{\infty}
\left( \delta \ln\frac{u}{\delta} \right)^n \right)\exp\left(-u \; \omega(N)\delta \left(1 + O\left(\delta\right)\right)
\sum_{n=1}^{\infty} n \left( \delta \ln\frac{u}{\delta} \right)^{n-1}\right)\right\}\, du.
\end{align*}
Next, we expand the exponentials and get (since $e^x =1 + x + O(x^{2})$ as $x \rightarrow 0$)
\begin{equation*}
K_{1}(\delta)=\int_{\delta}^{1/\sqrt{\delta}}\left[1-\left(j-2\right)\delta \ln\frac{u}{\delta}+\omega(N)\delta+O\left(\delta^{2} \ln^{2}\frac{u}{\delta}\right)\right]
\end{equation*}
\begin{equation*}
\qquad\qquad\quad\times e^{-u}\;\left\{1-u\sum_{n=1}^{\infty}\left(\delta \ln \frac{u}{\delta}\right)^{n}+O\left(u\sum_{n=1}^{\infty}\left(\delta \ln \frac{u}{\delta}\right)^{n}\right)^{2}\right\}
\end{equation*}
\begin{equation*}
\qquad\qquad\quad\times \left\{1-\omega(N)\;u\;\delta\left(1+O\left(\delta\right)\right)\sum_{n=1}^{\infty}\left(\delta \ln \frac{u}{\delta}\right)^{n-1}\right.
\end{equation*}
\begin{equation*}
\left.\qquad\qquad\qquad+O\left(\omega(N)\;u\;\delta\left(1+O\left(\delta\right)\right)\sum_{n=1}^{\infty}\left(\delta \ln \frac{u}{\delta}\right)^{n-1}\right)^{2}\right\}{du}
\end{equation*}
Hence,
\begin{equation*}
K_{1}(\delta)=\int_{\delta}^{1/\sqrt{\delta}}e^{-u}\left[1-\left(j-2+u\right)\delta \ln \frac{u}{\delta}+ \omega(N)\,\delta\left(1-u\right)+u^{2}O\left(\delta^{2}\ln^{2}\frac{u}{\delta}\right)\right]du
\end{equation*}
\begin{equation*}
=\int_{\delta}^{\infty}e^{-u}\left[1-\left(j-2+u\right)\delta \ln \frac{u}{\delta}+ \omega(N)\,\delta\left(1-u\right)+u^{2}O\left(\delta^{2}\ln^{2}\frac{u}{\delta}\right)\right]du
\end{equation*}
\begin{equation}
-\int_{1/\sqrt{\delta}}^{\infty}e^{-u}\left[1-\left(j-2+u\right)\delta \ln \frac{u}{\delta}+ \omega(N)\,\delta\left(1-u\right)+u^{2}O\left(\delta^{2}\ln^{2}\frac{u}{\delta}\right)\right]du.
\end{equation}\label{020}
However,
\begin{equation*}
\int_{1/\sqrt{\delta}}^{\infty}e^{-u}\left[1-\left(j-2+u\right)\delta \ln \frac{u}{\delta}+ \omega(N)\,\delta\left(1-u\right)+u^{2}O\left(\delta^{2}\ln^{2}\frac{u}{\delta}\right)\right]du
\end{equation*}
\begin{equation*}
=\int_{1/\sqrt{\delta}}^{\infty}
e^{-u}\left[1-u\left(\delta \ln\frac{u}{\delta}+\omega(N)\,\delta\right)+\omega(N)\,\delta-\left(j-2\right)\delta\,\ln\frac{u}{\delta}\right]{du}
\end{equation*}
\begin{equation*}
\leq\int_{1/\sqrt{\delta}}^{\infty}
e^{-u}\left[1-\left(1/\sqrt{\delta}\right)\left(\delta \ln\left(\frac{1/\sqrt{\delta}}{\delta}\right)+\omega(N)\,\delta\right)+\omega(N)\,\delta-\left(j-2\right)\delta\,\ln\left(\frac{1/\sqrt{\delta}}{\delta}\right)\right]{du}
\end{equation*}
\begin{equation}
\qquad\qquad\qquad+\int_{1/\sqrt{\delta}}^{\infty} u^{2} e^{-u}O\left(\delta^{2} \ln^{2}\frac{u}{\delta}\right){du}= O\left(e^{-1/\sqrt{\delta}}\right) \qquad \text{as\,\, $\delta\rightarrow 0^{+}$.}
\label{O}
\end{equation}
It follows that in the expression for $K_1(\delta)$ we can replace the upper limit of the integral by $\infty$. Therefore,  (\ref{eint5}) becomes (as $\delta\rightarrow 0^{+}$)
\begin{equation*}
I_{N}^{12} (\alpha;j )=\frac{\delta}{ f(N)^{j}}\int_{\delta}^{\infty}e^{-u}\left[1-\left(j-2+u\right)\delta \ln \frac{u}{\delta}+ \omega(N)\,\delta\left(1-u\right)+u^{2}O\left(\delta^{2}\ln^{2}\frac{u}{\delta}\right)\right]du
\end{equation*}
\begin{align}
=&\frac{\delta}{ f(N)^{j}}\left\{\left[1+\left(j-2\right)\delta \ln {\delta}+\omega(N)\,\delta\right]e^{-\delta}+\delta\left[ \ln{\delta}-\omega(N)\right]\left(1+\delta\right)e^{-\delta}\right\}\nonumber\\
-&\frac{\delta}{ f(N)^{j}}\left\{\delta \int_{\delta}^{\infty}e^{-u}\,u\,\ln u\,du+\left(j-2\right)\delta\int_{\delta}^{\infty}e^{-u}\,\ln{u}\,du\right\}\nonumber\\
+&\frac{\delta}{ f(N)^{j}}\,O\left(\int_{\delta}^{\infty}e^{-u}\,u^{2}\,\delta^{2}\ln^{2}\frac{u}{\delta}\,du\right).\label{25}
\end{align}
To continue we need some lemmas.\\\\
\textbf{Lemma 2.}
For the integral,
\begin{equation*}
M(x):=\int_{x}^{\infty}e^{-t}\,t\,\ln t\,dt,
\end{equation*}
we have the asymptotic expansion, as $x\rightarrow 0^{+}$,
\begin{equation}
M(x)\sim 1-\gamma+\frac{x^{2}}{2}\ln x-\frac{x^{2}}{4}+\cdots \label{M}
\end{equation}

{\it Proof}. Since
\begin{equation*}
\frac{dM(x)}{dx}=-x\ln x \;e^{-x}=-x\ln x\left(1-x+\frac{1}{2}x^{2}-\frac{1}{6}x^{3}+\cdots\right),
\end{equation*}
we have
\begin{equation}
M(x) \sim C_1 - \frac{1}{2} x^2 \ln x + \frac{1}{4}x^2 - \frac{1}{3} x^3 \ln x
- \frac{1}{9} x^3 - \frac{1}{8} x^4 \ln x + \frac{1}{32} x^4 + \cdots,
\label{G1a}
\end{equation}
where $C_1$ is a constant. In fact,
\begin{equation*}
C_1 = \int_{0}^{\infty}t\,\ln t \;e^{-t}dt.
\end{equation*}
Integration by parts yields
\begin{equation*}
C_1 = \int_0^{\infty} e^{-t}\,\ln t\,dt+\int_0^{\infty}  e^{-t} dt = \Gamma'(1)+1 =1 -\gamma
\end{equation*}
(see \cite{V}) and the proof is completed. \hfill $\blacksquare$ \\\\
\textbf{Lemma 3.}
For the integral,
\begin{equation*}
G(x):=\int_{x}^{\infty}\ln t \;e^{-t}dt,
\end{equation*}
we have the asymptotic expansion, as $x \rightarrow 0^{+}$,
\begin{equation}
G(x)\sim -\gamma-x\ln x+x+\frac{1}{2}x^{2}\ln x-\frac{1}{2}x^{2}-\frac{1}{6}x^{3}\ln x+\frac{1}{6}x^{3}+\frac{1}{24}x^{4}\ln x-\frac{1}{24}x^{4}+\cdots \label{G}
\end{equation}
The proof is similar to Lemma 2; it has been given in \cite{DP}.\\
We also observe that as $x\rightarrow 0^{+}$,
\begin {equation*}
\int_{x}^{\infty}t^{2}\; e^{-t} \ln^{2}t\; {dt} \sim C_0 :=
\int_{0}^{\infty}t^{2}\; e^{-t} \ln^{2}t\; {dt}=\frac{\pi^{2}}{3}+2\gamma^{2}-6\gamma+2,
\end{equation*}
and
\begin {equation*}
\int_{x}^{\infty}t^{2}\; e^{-t} \ln t\; {dt}\sim \tilde{C_{0}}:=
\int_{0}^{\infty}t^{2}\; e^{-t} \ln t\; {dt}=3-2\gamma.
\end{equation*}
In particular,
\begin {equation}
O\left(\int_{\delta}^{\infty}e^{-u}\,u^{2}\,\delta^{2}\ln^{2}\frac{u}{\delta}\,du\right)
= O\left(\delta^{2}\ln^{2}\delta\right)
\qquad
\text{as }\; \delta \rightarrow 0^+.
\label{127}
\end{equation}
Applying Lemmas 2, 3, and (\ref{127}), in (\ref{25}) we get (since $e^{-\delta}=1-\delta+O(\delta^{2})$ as $\delta \rightarrow 0^{+}$),
\begin{equation}
I_{N}^{12} (\alpha;j )=\frac{\delta}{ f(N)^{j}}\left[1+\left(j-1\right) \delta \ln \delta +\left[\left(j-1\right)\gamma-2  \right]\delta +O\left(\delta^{2} \ln^{2}\delta\right)\right].
\label{eint11}
\end{equation}
Notice that the error term in (\ref{eint11}) dominates the terms of (\ref{eint6}) and (\ref{O}).\\\\
\textbf{Remark 5.} In view of (\ref{delta}), (\ref{b3}) yields
\begin{equation}
I_{N}^{11} (\alpha;j ) << \frac{\delta^{3}\ln^{2}\delta}{f(N)^{j}},
\qquad j \geq 2,
\label{-10}
\end{equation}
as $N \to \infty$. Using (\ref{-10}) and (\ref{eint11}) and invoking (\ref{-11}), one has
\begin{equation}
I_N^1 (\alpha;j )=\frac{\delta}{ f(N)^{j}}\left[1+\left(j-1\right) \delta \ln \delta +\left[\left(j-1\right)\gamma-2  \right]\delta +O\left(\delta^{2} \ln^{2}\delta\right)\right].\label{000}
\end{equation}
\subsection{The integral $I_N^2 (\alpha;j )$}
Our next task is to compute the asymptotic behavior of the quantity $I_{N}^{2} (\alpha;j )$
defined in (\ref{I2}). It has been established in \cite{BP} that,
\begin{equation*}
\lim_{N}\sum_{i=1}^{N}\ln \left( 1-e^{-\frac{F(N)}{f(i)}s}\right) =0,
\end{equation*}
uniformly in $s\in [1,\infty)$.
From (\ref{I2}) we have
\begin{equation}
I_{N}^{2} (\alpha;j )= \sum_{k=1}^{N}f(k)^{-j}\int_1^{\infty}s^{j-1}\frac{e^{-\frac{F(N)}{f(k)}s}}{1-e^{-\frac{F(N)}{f(k)}s}}\;\exp \left[ \sum_{i=1}^{N}\ln \left( 1-e^{-\frac{F(N)}{f(i)}s}\right) \right] ds.\label{001}
\end{equation}
Using the Taylor expansion of the logarithm, the comparison of sums and integrals (i.e. (\ref{0})), and Lemma 1, for $m=0,$ we get
\begin{equation*}
I_{N}^{2} (\alpha;j )= \sum_{k=1}^{N}f(k)^{-j}\int_1^{\infty}s^{j-1}\frac{e^{-\frac{F(N)}{f(k)}s}}{1-e^{-\frac{F(N)}{f(k)}s}}
\left[1+O\left(\frac{f(N)^{2}}{F(N)f^{\prime}(N)}e^{-\frac{F(N)}{f(N)}s}\right)\right]ds.
\end{equation*}
The scaling $u=F(N)f(k)^{-1}s,$ via the definition of $F$ yields
\begin{align}
\int_1^{\infty}s^{j-1}\frac{e^{-\frac{F(N)}{f(k)}s}}{1-e^{-\frac{F(N)}{f(k)}s}}ds&=
\left(\frac{f(k)}{F(N)}\right)^{j}\int_{\frac{F(N)}{f(k)}}^{\infty}\frac{e^{-u}}{1-e^{-u}}\;u^{j-1}\;du\nonumber\\
&=\frac{f(k)}{F(N)}e^{-\frac{F(N)}{f(k)}}+O\left[\left(\frac{f(k)}{F(N)}\right)^{2}e^{-\frac{F(N)}{f(k)}}\right].\label{26}
\end{align}
Thus, $I_{N}^{2} (\alpha;j )$ (see (\ref{001})), becomes
\begin{equation*}
I_{N}^{2} (\alpha;j )= \sum_{k=1}^{N}f(k)^{-j}\left(\frac{f(k)}{F(N)}e^{-\frac{F(N)}{f(k)}}+O\left[\left(\frac{f(k)}{F(N)}\right)^{2}e^{-\frac{F(N)}{f(k)}}\right]\right)
\end{equation*}
\begin{equation*}
=\frac{1}{F(N)}\sum_{k=1}^{N}\frac{1}{f(k)^{j-1}}e^{-\frac{F(N)}{f(k)}}+O\left[\frac{1}{F(N)^{2}}\sum_{k=1}^{N}\frac{1}{f(k)^{j-2}}e^{-\frac{F(N)}{f(k)}}\right].
\end{equation*}
Under conditions (\ref{C1}), for sufficiently large  $N$ we have
\begin{equation*}
\sum_{k=1}^{N}\frac{1}{f(k)^{j-1}}e^{-\frac{F(N)}{f(k)}}=\int_{1}^{N}\frac{1}{f(x)^{j-1}}e^{-\frac{F(N)}{f(x)}}\,dx+O\left(\frac{1}{f(N+1)^{j-1}}e^{-\frac{F(N)}{f(N+1)}}\right).
\end{equation*}
Applying Lemma 1, for $m=(1-j),$ $s=1,$ yields
\begin{align*}
\sum_{k=1}^{N}\frac{1}{f(k)^{j-1}}e^{-\frac{F(N)}{f(k)}}&=
\frac{f(N)^{3-j}}{F(N)f^{\prime}(N)}e^{-\frac{F(N)}{f(N)}}+\omega(N)\;\frac{f(N)^{4-j}}{F(N)^{2}f^{\prime}(N)}e^{-\frac{F(N)}{f(N)}}\left[1+O\left(\frac{f(N)}{F(N)}\right)\right]\nonumber\\
&+O\left(\frac{1}{f(N+1)^{j-1}}e^{-\frac{F(N)}{f(N+1)}}\right).
\end{align*}
By (\ref{LIMIT}) the above quantity becomes
\begin{equation*}
\sum_{k=1}^{N}\frac{1}{f(k)^{j-1}}e^{-\frac{F(N)}{f(k)}}=\frac{f(N)^{3-j}}{F(N)f^{\prime}(N)}e^{-\frac{F(N)}{f(N)}s}+\omega(N)\;\frac{f(N)^{4-j}}{F(N)^{2}f^{\prime}(N)}e^{-\frac{F(N)}{f(N)}s}\left[1+O\left(\frac{f(N)}{F(N)}\right)\right].
\end{equation*}
Hence,
\begin{equation*}
I_{N}^{2} (\alpha;j )=\frac{f(N)^{3-j}}{F(N)^{2}f^{\prime}(N)}e^{-\frac{F(N)}{f(N)}}\left[1+O\left(\frac{f(N)}{F(N)}\right)\right].
\end{equation*}
Using the definition of $F$, i.e. (\ref{F}), and (\ref{delta}) one has
\begin{equation}
I_{N}^{2} (\alpha; j) =\frac{1}{f(N)^{j}\ln^{j}\left(\frac{f(N)}{f^{\prime}(N)}\right)}\left[1+O\left(\frac{1}{\ln\frac{f(N)}{f^{\prime}(N)}}\right)\right]=\frac{\delta}{f(N)^{j}}\left(\delta^{j-1}+O\left(\delta^{j}\right)
\right).
\label{27}
\end{equation}
We are now ready for our main result.

\subsection{Conclusion. Asymptotics of $E[U_j^N]$}
Recall that when the main collector has completed her album, they remain $U_{j}^{N}$ unfilled places in the album of the
$j$-th collector, $j = 2, 3, \dots$. The asymptotics of $E\,[\,U_{j}^{N}\,]$ is given by the following theorem.\\\\
\textbf{Theorem 1.}
Let $\alpha = \{a_k\}_{k=1}^{\infty } = \{1/f(k)\}_{k=1}^{\infty}$, where $f$ satisfies (\ref{C1a}) and (\ref{C1}). If the coupon frobabilities (i.e. the $p_k$'s) are as in (\ref{B0}), then, as $N \rightarrow \infty$, we have
\begin{equation}
E\left[\,U_{2}^{N}\,\right]=\frac{1}{\delta}\left[1+ \delta \ln \delta +\left(\gamma-1  \right)\delta +O\left(\delta^{2} \ln^{2}\delta\right)\right],\label{FINAL1}
\end{equation}
\begin{equation}
E\left[\,U_{j}^{N}\,\right]=\frac{1}{\left(j-1\right)!\,\,\delta^{j-1}}\left[1+\left(j-1\right) \delta \ln \delta +\left[\left(j-1\right)\gamma-2  \right]\delta + O\left(\delta^{2} \ln^{2}\delta\right)\right]
\label{FINAL2}
\end{equation}
for $j \geq 3$. Recall that
\begin{equation*}
\delta = \frac{1}{\ln \left(f(N) / f^{\prime}(N)\right)}.
\end{equation*}
{\it Proof}. The desired result follows by using formulas (\ref{000}) and (\ref{27}) in (\ref{20}) and (\ref{5}). Notice that, for $j \geq 3$ all three terms of the asymptotics in (\ref{FINAL2}) come solely from
$I_N^{12} (\alpha; j)$ (see (\ref{eint11})). For $j = 2$ part of the third term
of the asymptotics in (\ref{FINAL1}) is due to $I_N^2 (\alpha; j)$ (see (\ref{27})), while the rest of the third term, as well as the first two terms
are, again, due to $I_N^{12} (\alpha; j)$ (see (\ref{eint11})). The integral
$I_N^{12} (\alpha; j)$ does not contribute at all in the first three terms of
the asymptotics of $E[U_j^N]$.
\hfill $\blacksquare$\\\\
\textbf{Remark 6.} From Theorem 1 and for all $j \geq 2$ we have
\begin{equation}
\,\,\,\,\,E\left[U_{j}^{N}\right]\sim\frac{1}{\left(j-1\right)!}\ln\left(\frac{f(N)}{f^{\prime}(N)}\right)^{j-1},
\,\,\,\,\,\, N\rightarrow \infty.
\label{FINAL3}
\end{equation}

\noindent \textbf{Example 1.} $a_k = 1/k^p$, where $p > 0$. This is the so-called \textit{generalized Zipf} law (for detailed asymptotic results regarding the first collector, i.e. the random variable $T_N$, see \cite{DP}). These decaying sequences fall, clearly, into the previous discussion, since $f(x) = x^p$, satisfies
(i)--(iv) of (\ref{C1}). Here,
\begin{equation}
\delta = \frac{1}{\ln N - \ln p}
\label{DDD1}
\end{equation}
and, hence, Theorem 1 gives
\begin{equation}
E\left[U_2^N\right] = \ln N - \ln(\ln N) + (\gamma - 1 - \ln p) + O\left(\frac{\ln(\ln N)^2}{\ln N}\right)
\label{30}
\end{equation}
and
\begin{align}
E\left[U_j^N\right]&
= \frac{(\ln N)^{j-1}}{(j-1)!} - \frac{(\ln N)^{j-2} \ln(\ln N)}{(j-2)!} + \frac{(j - 1)(\gamma - \ln p) - 2}{(j-1)!}(\ln N)^{j-2}
\nonumber\\
& \qquad \qquad \qquad \qquad \quad \,\,\,\,
+ O\left(\left(\ln N \right)^{j-3}\ln (\ln N)^2\right)
\label{DDD30}
\end{align}
for $j \geq 3$. Notice that $p$ does not appear in the first two terms of the asymptotics. Also, the leading term
is the same as in the case where all the $p_k$'s are equal (see (\ref{15})), however the second terms of the asymptotics differ.\\\\
\noindent \textbf{Example 2.} $a_k = \exp(-p k^q)$, where $p > 0$ and
$0 < q < 1$. Again, these decaying sequences fall into
the previous discussion, since $f(x) = \exp(p x^q)$, satisfies (i)--(iv) of (\ref{C1}). Here,
\begin{equation*}
\delta = \frac{1}{\ln(N^{1-q}) - \ln(pq)}.
\end{equation*}
If we compare the above formula with (\ref{DDD1}) we can see that the
asymptotics of $E[U_j^N]$ ($j \geq 2$) for Example 2 can be obtained from the
formulas (\ref{30}) and (\ref{DDD30}) of Example 1 after replacing $N$ by
$N^{1-q}$ and $p$ by $pq$. For example (this, also follows from (\ref{FINAL3})),
\begin{equation*}
E\left[U_j^N\right] \sim \frac{(\ln(N^{1-q})^{j-1}}{(j-1)!}
= (1 - q)^{j-1} \frac{(\ln N)^{j-1}}{(j-1)!},
\qquad
N \to \infty,
\end{equation*}
for all $j \geq 2$. Notice that $p$ does not appear in the leading asymptotics of $E[U_j^N]$; it first appears in the third term.

\section{Growing sequences}
In this section we will examine sequences
$\alpha = \{a_k\}_{k=1}^{\infty}$, such that
\begin{equation*}
a_k \rightarrow \infty.
\end{equation*}
We will exhibit several cases where (recall (\ref{5})),
$E[U_j^N] = I_N(\alpha; j)$ approaches a finite limit as $N \rightarrow \infty$.
However, we will also see that there are cases for which
$\lim_N E[U_j^N] = \infty$.

By substituting $x = e^{-t}$, (\ref{4}) becomes
\begin{equation}
I_N(\alpha; j) =
\frac{1}{(j - 1)!}
\int_0^1 \left(\sum_{k=1}^N a^j_k \, \frac{x^{a_k}}{1 - x^{a_k}}\right)
\,\left[\prod_{k=1}^N \left(1 - x^{a_k}\right)\right]
\,|\ln x|^{j-1}\,\frac{dx}{x}.
\label{4b}
\end{equation}
Fot typographical convenience let us set
\begin{equation}
S_N(x) = S_N(x; \alpha) := \sum_{k=1}^N x^{a_k},
\label{G1}
\end{equation}
\begin{equation}
F_N(x) = F_N(x; \alpha) := \prod_{k=1}^N \left(1 - x^{a_k}\right),
\label{G2}
\end{equation}
and
\begin{equation}
L_N(x; j) = L_N(x ;\alpha; j) :=
\sum_{k=1}^N a^j_k\,\frac{x^{a_k}}{1 - x^{a_k}},
\qquad
j = 2, \dots, r
\label{G3}
\end{equation}
(in fact, $L_N(x ;\alpha; j$) makes sense for any real number $j$). Obviously, for a fixed $x \in (0, 1)$ and a fixed $j$ we have that $S_N(x)$ and $L_N(x; j)$ increase with $N$, while $F_N(x)$ decreases.\\
Given a sequence $\alpha = \{a_k\}_{k=1}^{\infty}$ of positive terms let
\begin{equation}
x_{\alpha} := \inf\left\{x \in [0, 1] \, : \,
\sum_{k=1}^{\infty} x^{a_k} = \infty\right\}.
\label{GS1}
\end{equation}
From now on, we will consider only sequences $\alpha$ such that
\begin{equation}
0 < x_{\alpha} \leq 1.
\label{GS2}
\end{equation}
Roughly speaking, condition (\ref{GS1})--(\ref{GS2}) says that $a_k$ grows at least logarithmically. For example, if $a_k = \ln k$ ($k \geq 2$), then
$x_{\alpha} = 1/e$ and (\ref{GS2}) is satisfied. However, if
$a_k = \epsilon_k \ln k$, where $\epsilon_k \to 0$, then $x_{\alpha} = 0$, i.e. (\ref{GS2}) is not satisfied. To have $x_{\alpha} = 1$, $a_k$ must grow faster than $\ln k$ (roughly speaking). For instance, if $a_k = \lambda_k \ln k$, where $\lambda_k \to \infty$, then $x_{\alpha} = 1$.\\
Let us set
\begin{equation}
S(x) = S(x; \alpha) := \sum_{k=1}^{\infty} x^{a_k}.
\label{GS1a}
\end{equation}
Then, from condition (\ref{GS1})--(\ref{GS2}) it follows that $S(x) < \infty$ for
$x \in [0, x_{\alpha})$ and $S(x) = \infty$ for $x \in (x_{\alpha}, 1]$. Furthermore, if $x_{\alpha} < 1$, then, depending on $\alpha$, $S(x_{\alpha})$ can be finite or infinite. For instance, if $a_k = \ln k$ ($k \geq 2$), then $x_{\alpha} = 1/e$ and $S(1/e) = \infty$, while if $a_k = \ln(k \ln^2 k)$ ($k \geq 3$), then, again $x_{\alpha} = 1/e$, but now $S(1/e)$ is finite. Of course,
if $x_{\alpha} = 1$, then $S(x_{\alpha}) = S(1) = \infty$.\\
If we set
\begin{equation}
F(x) = F(x; \alpha) := \lim_N F_N(x; \alpha)
= \prod_{k=1}^{\infty} \left(1 - x^{a_k}\right),
\label{G5}
\end{equation}
then, in view of (\ref{GS1})--(\ref{GS2}), by standard properties of infinite products (see, e.g., \cite{Ru}) we have that
\begin{equation}
0 < F(x) < 1
\qquad
\text{for all }\; x \in (0, x_{\alpha})
\label{GS3}
\end{equation}
and, also, $F(x) = 0$ for $x \in (x_{\alpha}, 1]$.
Of course, $F(0) = 1$ and $F(1) = 0$. Furthermore, $F(x)$ is (decreasing on
$[0, 1]$ and) continuous for all $x \in [0, 1]$ with only one possible exception (at $x= x_{\alpha}$): $F(x)$ is not continuous at $x_{\alpha}$ if and only if $S(x_{\alpha}) < \infty$, since in this case
$F(x_{\alpha}-) = F(x_{\alpha}) > 0$, while $F(x_{\alpha}+) = 0$. Of course, if $x_{\alpha} = 1$, then $F(x)$ is continuous on $[0, 1]$.\\
Next, we notice that it is easy to show that, under (\ref{GS1})--(\ref{GS2}) we have
\begin{equation}
L(x; j) = L(x; \alpha; j) := \lim_N L_N(x; \alpha; j)
= \sum_{k=1}^{\infty} a^j_k\,\frac{x^{a_k}}{1 - x^{a_k}} < \infty,
\quad
x \in (0, x_{\alpha}).
\label{G7}
\end{equation}
Furthermore,
\begin{equation}
L(x; j) = O(x^{a_{\text{min}}})
\qquad
\text{as }\; x \to 0^+,
\label{G8}
\end{equation}
where
\begin{equation}
a_{\min} := \min_{k \in \mathbb{N}} (a_k) > 0
\label{G8b}
\end{equation}
(the strict positivity follows from the fact that
$a_k > 0$ for all $k \in \mathbb{N}$ and $a_k \to \infty$).
As for the value of $L(x_{\alpha}; j)$, depending on the sequence $\alpha$ there might be a $j_0$ such that $L(x_{\alpha}; j) < \infty$, for $j < j_0$, while $L(x_{\alpha}; j) = \infty$, for $j \geq j_0$. For example, if $a_k = \ln(k \ln^p k)$ ($k \geq 3$), then $x_{\alpha} = 1/e$ and $L(1/e; j) < \infty$ for $j < p-1$, while $L(1/e; j) = \infty$ for $j \geq p-1$.\\
Inspired by (\ref{4b}) we introduce the quantity
\begin{align}
I(\alpha; j) :=& \frac{1}{(j - 1)!}
\int_0^{x_{\alpha}} L(x; \alpha; j) F(x; \alpha)\,|\ln x|^{j-1}\,\frac{dx}{x}
\label{GS4}
\\
=& \frac{1}{(j - 1)!}
\int_{-\ln(x_{\alpha})}^\infty
L\left(e^{-t}; \alpha; j\right) F\left(e^{-t}; \alpha\right) \, t^{j-1} dt.
\label{GS4b}
\end{align}

\noindent \textbf{Proposition 2.} If $I(\alpha; j) = \infty$ (see (\ref{GS4})),
then
\begin{equation*}
\lim_N E[U_j^N] = \infty.
\end{equation*}
\textit{Proof}. First notice that by (\ref{5}), (\ref{4b}), (\ref{G2}), and
(\ref{G3}) we have
\begin{align*}
E[U_j^N] = I_N(\alpha; j) &=
\frac{1}{(j - 1)!}
\int_0^1 L_N(x; \alpha; j) F_N(x; \alpha)\,|\ln x|^{j-1}\,\frac{dx}{x}
\\
&\geq
\frac{1}{(j - 1)!}
\int_0^{x_{\alpha}} L_N(x; \alpha; j) F_N(x; \alpha)\,|\ln x|^{j-1} \, \frac{dx}{x}\, .
\end{align*}
Also, by (\ref{G5}) and (\ref{G7})
\begin{equation*}
\lim_N \left[L_N(x; \alpha; j) F_N(x; \alpha)\right]
= L(x; \alpha; j) F(x; \alpha),
\qquad
0 \leq x < x_{\alpha}.
\end{equation*}
Thus, the proposition follows from (\ref{GS4}) and the Fatou Lemma.
\hfill $\blacksquare$\\

\subsection{Some results for the case $x_{\alpha} = 1$}

\noindent \textbf{Proposition 3.} Assume that for the sequence
$\alpha = \{a_k\}_{k=1}^{\infty}$ we have $x_{\alpha} = 1$ (recall (\ref{GS1})).
If there is an $N_0$ such that
\begin{equation}
\int_0^1 L(x; \alpha; j) F_{N_0}(x; \alpha)\,|\ln x|^{j-1}\,\frac{dx}{x}
< \infty,
\label{G10}
\end{equation}
then (recalling (\ref{GS4}))
\begin{equation}
\lim_N E[U_j^N] = I(\alpha; j) < \infty.
\label{G11}
\end{equation}
\textit{Proof}. As in the proof of Proposition 2
\begin{equation*}
E[U_j^N] = I_N(\alpha; j) =
\frac{1}{(j - 1)!}
\int_0^1 L_N(x; \alpha; j) F_N(x; \alpha)\,|\ln x|^{j-1}\,\frac{dx}{x}.
\end{equation*}
and
\begin{equation*}
\lim_N \left[L_N(x; \alpha; j) F_N(x; \alpha)\right]
= L(x; \alpha; j) F(x; \alpha),
\qquad
0 \leq x < 1.
\end{equation*}
Since
\begin{equation*}
0 < L_N(x; \alpha; j) F_N(x; \alpha) < L(x; \alpha; j) F_{N_0}(x; \alpha)
\qquad \text{for }\; N > N_0.
\end{equation*}
(\ref{G11}) follows by dominated convergence.
\hfill $\blacksquare$\\\\
If $\alpha = \{a_k\}_{k=1}^{\infty}$ (with $x_{\alpha} = 1$) satisfies
$a_k = o(k^p)$ for all $p > 0$  then condition (\ref{G10}) cannot be satisfied.
As an example of such a sequence one can take $a_k = (\ln k)^q$, $q > 1$. In this case Proposition 3 is inconclusive.\\

\noindent \textbf{Remark 7.} Suppose that for the sequence $\alpha = \{a_k\}_{k=1}^{\infty}$ we have $x_{\alpha} = 1$. Let
$\beta = \{b_k\}_{k=1}^{\infty}$ be a sequence such that there is an integer $k_0 \geq 0$ for which
\begin{equation}
b_k = a_{k+k_0},
\qquad
k \in \mathbb{N}
\label{G12}
\end{equation}
(i.e. $\beta$ is the $k_0$-left shift of $\alpha$). Then $x_{\beta} = 1$ and it is easy to see that $\alpha$ satisfies condition (\ref{G10}) for some $N_0 = N_0(\alpha)$ if and only if $\beta$ satisfies condition (\ref{G10}) for some $N_0 = N_0(\beta)$.\\
The same equivalence is true if given $\alpha$ the sequence
$\tilde{\beta} = \{\tilde{b}_k\}_{k=1}^{\infty}$ is such that
\begin{equation}
|a_k - \tilde{b}_k| \leq M,
\qquad
k \in \mathbb{N},
\label{G13}
\end{equation}
for some $M > 0$.
In particular, for a given sequence $\alpha = \{a_k\}_{k=1}^{\infty}$, if we set $\underline{\alpha} :=
\{\left\lfloor a_k \right\rfloor\}_{k=1}^{\infty}$ and
$\overline{\alpha} :=
\{\left\lceil  a_k \right\rceil\}_{k=1}^{\infty}$ (where
$\left\lfloor x \right\rfloor$ and $\left\lceil  x \right\rceil$ denote the greatest integer $\leq x$ and the smallest integer $\geq x$ respectively), then
$\alpha$ satisfies condition (\ref{G10}) if and only if $\underline{\alpha}$ satisfies condition (\ref{G10}) (and this in turn holds if and only if $\overline{\alpha}$ satisfies condition (\ref{G10})).\\\\
Let us now assume that all the terms of the sequence
$\alpha = \{a_k\}_{k=1}^{\infty}$ are positive \emph{integers}. Then $L(x; \alpha; j)$ of (\ref{G7}) can be expressed as
\begin{equation}
L(x; \alpha; j) = \sum_{m=1}^{\infty} m^j A(m)\frac{x^m}{1 - x^m}
\label{G14}
\end{equation}
where
\begin{equation}
A(m) := \#\left\{a_k \,:\, a_k = m \right\}
\label{G14a}
\end{equation}
(the symbol $\#$ indicates cardinality). The series in the right-hand side of (\ref{G14}) is a so-called Lambert series (see, e.g., \cite{HW}) and can be easily transformed to a power series. Indeed,
\begin{equation}
L(x; \alpha; j) = \sum_{n=1}^{\infty} A_L(n) x^n,
\label{G15}
\end{equation}
where
\begin{equation}
A_L(n) := \sum_{d|n} d^j A(d),
\label{G16}
\end{equation}
i.e. the sum is taken over all divisors $d$ of $n$.\\

\noindent \textbf{Corollary 1.} Suppose that the sequence $\alpha = \{a_k\}_{k=1}^{\infty}$ has integer terms and
satisfies (recall (\ref{G7}), (\ref{G14a}), (\ref{G15}), and (\ref{G16}))
\begin{equation}
L(x; \alpha; j) = \sum_{n=1}^{\infty} A_L(n) x^n = O\left( \frac{1}{(1 - x)^{\rho}} \right),
\qquad
x \to 1^-,
\label{G17}
\end{equation}
for some $\rho > 0$. Then $\alpha$ satisfies (\ref{G10}) and, consequently,
\begin{equation*}
\lim_N E[U_j^N] = I(\alpha; j) < \infty.
\end{equation*}
\textit{Proof}. If (\ref{G17}) is true, then by choosing $N_0 = \left\lceil \rho \right\rceil$ we get that (see (\ref{G2}))
\begin{equation*}
L(x; \alpha; j) F_{N_0}(x; \alpha)
\qquad
\text{is bounded for }\; x \in [0, 1),
\end{equation*}
hence (\ref{G10}) is satisfied (with the help of (\ref{G8}), which takes care of the lower limit, i.e. $x = 0$, of the integral in (\ref{G10})).
\hfill $\blacksquare$\\

\noindent \textbf{Corollary 2.} Let $\alpha = \{a_k\}_{k=1}^{\infty}$ be a sequence of positive integers for which $A(m)$ of (\ref{G14a})) satisfies
\begin{equation}
A(m) = O(m^\nu)
\label{G18}
\end{equation}
for some $\nu > 0$. Then
\begin{equation*}
\lim_N E[U_j^N] = I(\alpha; j) < \infty
\qquad \text{for all }\; j \geq 2.
\end{equation*}
\textit{Proof}. From (\ref{G16}) and (\ref{G18}) we get that
\begin{equation}
A_L(n) = \sum_{d|n} d^j A(d) \leq \sum_{d=1}^n d^j A(d) = O(n^{\nu +j+1}),
\label{G18a}
\end{equation}
Set $\mu = \left\lceil \nu +j+1 \right\rceil$. Then, for $x \in [0, 1)$ condition (\ref{G18a}) implies
\begin{equation*}
L(x; \alpha; j) \leq M \sum_{n=0}^{\infty} (n+1)(n+2) \cdots (n+\mu) x^n
= M \frac{d^{\mu}}{dx^{\mu}} \left[ \frac{1}{1 - x} \right]
\end{equation*}
for some constant $M > 0$. In other words
\begin{equation*}
L(x; \alpha; j) \leq \frac{M \mu !}{(1 - x)^{\mu +1}}
\end{equation*}
and hence (\ref{G17}) of Corollary 1 is satisfied by choosing $\rho = \mu +1$.
\hfill $\blacksquare$\\

\noindent \textbf{Remark 8.} Suppose that $\alpha = \{a_k\}_{k=1}^{\infty}$ is a sequence whose terms are positive reals, not necessarily integers. Let
\begin{equation}
A^*(m) := \#\left\{a_k \,:\, a_k \leq m \right\}.
\label{G16b}
\end{equation}
Notice that
\begin{equation}
A^*(m) \geq
\#\left\{a_k \,:\, \left\lceil a_k \right\rceil \leq m \right\}.
\label{G16c}
\end{equation}
If  there is a $\nu > 0$ such that
\begin{equation}
A^*(m) = O(m^\nu)
\label{G16d}
\end{equation}
then (\ref{G16c}) implies that
\begin{equation*}
\#\left\{a_k \,:\, \left\lceil a_k \right\rceil \leq m \right\}  = O(m^\nu),
\end{equation*}
i.e. $\overline{\alpha} =
\{\left\lceil a_k \right\rceil\}_{k=1}^{\infty}$ satisfies (\ref{G18}).
Therefore, with the help of the last part of Remark 7 we deduce that $\lim_N E[U_j^N] = I(\alpha; j) < \infty$ for all $j \geq 2$.\\

\noindent \textbf{Example 3.} Let $a_k = k^p$, where $p > 0$ (the case $p = 1$ is known as the \textit{linear} case). Then (recall (\ref{G16b}))
\begin{equation*}
A^*(m) = \#\left\{k \in \mathbb{N}: k^p \leq m \right\} = O(m^{1/p}).
\end{equation*}
Thus, by Remark 8 we obtain that
\begin{equation*}
\lim_N E[U_j^N] = I(\alpha; j) < \infty,
\qquad \text{for all }\; j \geq 2.
\end{equation*}

\noindent \textbf{Example 4.} Let $a_k = e^{p k^q}$, where $p, q > 0$. Then,
we can, again, use Remark 8 as in Example 3 to conclude that
\begin{equation*}
\lim_N E[U_j^N] = I(\alpha; j) < \infty,
\qquad \text{for all }\; j \geq 2.
\end{equation*}
In the same way we can see that for the sequence $a_k = k!$, $k = 1, 2, \dots$,
we also have $\lim_N E[U_j^N] = I(\alpha; j) < \infty$ for all $j \geq 2$.\\
Let us, also, discuss the sequence $\beta = \{ b_k \}_{k=1}^{\infty}$,
with $b_k = e^{-pk}$, $p > 0$. In this case, the function $f(x) = e^{px}$ does not satisfy condition (ii) of (\ref{C1}), thus Theorem 1 cannot be applied.
However, the sequences $\alpha =  \{ e^{pk} \}_{k=1}^{\infty}$
and $\beta$ produce the same coupon
probabilities! This follows from the fact that for each $N$, if we let
$c_N = e^{pN}$, then $\left\{ a_{j}:\,0\leq j\leq N\right\} =
\left\{c_{N}b_{j}:\,0\leq j\leq N\right\}$, i.e. the elements of the two truncated sequences are proportional to each others. It follows that
\begin{equation*}
\lim_N E[U_j^N] = I(\alpha; j) < \infty,
\qquad \text{for all }\; j \geq 2.
\end{equation*}
Here, $I(\alpha; j)$ depends on $p$ (compare with Example 2).\\

\subsection{Two examples for the case $0 < x_{\alpha} < 1$}

\noindent \textbf{Example 5.} If $a_k = \ln(k \ln^2 k) = \ln k + 2 \ln(\ln k)$
($k \geq 3$), then
\begin{equation*}
\lim_N E[U_j^N] = \infty
\qquad \text{for all }\; j \geq 2.
\end{equation*}
In order to justify this equation let us observe that, by Proposition 2 it suffices to show that
$I(\alpha; j) = \infty$, where $I(\alpha; j)$ is defined in (\ref{GS4})--(\ref{GS4b}). First we notice that $x_{\alpha} = 1/e$ and $S(1/e) < \infty$.\\
Clearly, $1 > F\left(e^{-t}\right) \geq F(1/e) > 0$. Thus, from (\ref{GS4b}) we have
\begin{equation*}
I(\alpha; j) =
\frac{1}{(j-1)!} \int_1^{\infty} t^{j-1} \left(\sum_{k=3}^{\infty}
\frac{a^j_k\,e^{-a_k t}}{1 - e^{-a_k t}}\right)\,F\left(e^{-t}\right)\,dt.
\end{equation*}
Since $1 / (1 - e^{-a_k t}) > 1$ for all $t > 0$, we have (with the help of Tonelli's theorem)
\begin{align*}
I(\alpha; j)\geq \,\, &
\frac{F(1/e)}{(j-1)!}
\int_1^{\infty} t^{j-1}
\left(\sum_{k=3}^{\infty} a^j_k \,e^{-a_k t}\right)\,dt
\\
\geq \,\,&
\frac{F(1/e)}{(j-1)!}
\int_1^{\infty}
\left(\sum_{k=3}^{\infty} a^j_k \,e^{-a_k t}\right)\,dt
\\
=\,\, &
\frac{F(1/e)}{(j-1)!} \,
\sum_{k=3}^{\infty}
a^j_k \left(\int_1^{\infty} e^{-a_k t}\right)\,dt
\\
=\,\, &
\frac{F(1/e)}{(j-1)!} \, \sum_{k=3}^{\infty} a^{j-1}_k\, e^{-a_k},
\end{align*}
and substituting $a_k =\ln(k \ln^2 k)$ we get
\begin{equation*}
I(\alpha; j) \geq
\frac{F(1/e)}{(j-1)!}\,
\sum_{k=3}^{\infty} \frac{\left[\ln(k \ln^2 k)\right]^{j-1}}{k \ln^2 k}
= \infty,
\qquad j \geq 2.
\end{equation*}
Finally, it is instructive to compare Example 5 with the following example.\\

\noindent \textbf{Example 6.} Let $a_k = \ln k$ ($k\geq 3$). Then,
$x_{\alpha} = 1/e$ and $S(1/e) = \infty$ (equivalently,
$S\left(e^{-t}\right) = \sum_{k=2}^{\infty} e^{-t a_k} < \infty$ if and only if $t > 1$). Here, in contrast with Example 5, we will show that, for all
$j \geq 2$, $E[U_j^N]$ approaches a finite limit as $N \to \infty$.\\
From (\ref{5}) and (\ref{4}) we have
\begin{equation}
E[U_j^N]
= \frac{1}{(j-1)!} \left[Q_N^1(\alpha; j) + Q_N^{\infty}(\alpha; j)\right],
\label{H1}
\end{equation}
where
\begin{equation}
Q_N^1(\alpha; j) :=
\int_0^1 \left[\sum_{k=3}^N a_k^j \,
\frac{e^{-a_k t}}{1 - e^{-a_k t}} \, t^{j-1} \right]
\left[\prod_{k=3}^N \left(1 - e^{-a_k t}\right)\right] dt
\label{H2}
\end{equation}
and
\begin{equation}
Q_N^{\infty}(\alpha; j) :=
\int_1^{\infty} \left[\sum_{k=3}^N a_k^j \,
\frac{e^{-a_k t}}{1 - e^{-a_k t}} \, t^{j-1} \right]
\left[\prod_{k=3}^N \left(1 - e^{-a_k t}\right)\right] dt.
\label{H3}
\end{equation}
Notice that $a_k = \ln k$ yields $e^{-a_k t} = k^{-t}$. Let us analyze $Q_N^1(\alpha; j)$ first. Since
\begin{equation*}
\frac{x}{1 - e^{-x}} < 1 + x,
\qquad x \in (0, \infty),
\end{equation*}
formula (\ref{H2}) implies
\begin{equation*}
Q_N^1(\alpha; j) <
\int_0^1 \left[\sum_{k=3}^N a_k^{j-1} \,
e^{-a_k t} (1 + a_k t) \, t^{j-2} \right]
\left[\prod_{k=3}^N \left(1 - e^{-a_k t}\right)\right] dt
\end{equation*}
or (since $a_k = \ln k$)
\begin{align}
Q_N^1(\alpha; j) &<
\int_0^1 \left[\sum_{k=3}^N \left(\frac{(\ln k)^{j-1}}{k^t} \, t^{j-2}
+ \frac{(\ln k)^j}{k^t} \, t^{j-1} \right) \right]
\left[\prod_{k=3}^N \left(1 - \frac{1}{k^t}\right)\right] dt
\nonumber
\\
&< 2 \int_0^1 \left[\sum_{k=3}^N \frac{(\ln k)^j}{k^t}\right] \, t^{j-2} \,
\left[\prod_{k=3}^N \left(1 - \frac{1}{k^t}\right)\right] dt.
\label{H5}
\end{align}
The expression $(\ln k)^j / k^t$ (viewed as a function of $k$) has a unique maximum. It is attained when $\ln k = j/t$ and the maximum value is
$j^j e^{-j} / t^j$. It follows that
\begin{equation*}
\sum_{k=3}^N \frac{(\ln k)^j}{k^t}
< \frac{j^j e^{-j}}{t^j}
+ \int_3^N \frac{(\ln \kappa)^j}{\kappa^t} \, d\kappa
< \frac{j^j e^{-j}}{t^j}
+ \int_0^{\ln N} \xi^j e^{(1-t) \xi} d\xi.
\end{equation*}
The integrand in the second integral above is increasing in $\xi$. Hence the integral is bounded by the value of the integrand at $\xi = \ln N$ times $\ln N$ (i.e the length of the interval of integration). Thus
\begin{equation}
\sum_{k=3}^N \frac{(\ln k)^j}{k^t}
< \frac{j^j e^{-j}}{t^j}
+ (\ln N)^{j+1} N^{1-t}.
\label{H6}
\end{equation}
Using (\ref{H6}) in (\ref{H5}) we get
\begin{equation}
Q_N^1(\alpha; j)
< \int_0^1 \frac{1}{t^2} \, \left[c_j
+ (\ln N)^{j+1} N^{1-t} t^j \right]
\left[\prod_{k=3}^N \left(1 - \frac{1}{k^t}\right)\right] dt,
\label{H7}
\end{equation}
where $c_j := j^j e^{-j}$. Hence,
\begin{equation}
Q_N^1(\alpha; j)
< c_j R_N(\alpha; j) + Z_N(\alpha; j),
\label{H8}
\end{equation}
where
\begin{align}
R_N(\alpha; j) :=&
\int_0^1 \frac{1}{t^2} \,
\left[\prod_{k=3}^N \left(1 - \frac{1}{k^t}\right)\right] dt
\nonumber
\\
=& \int_0^1 \frac{(1 - e^{-t \ln 3}) (1 - e^{-t \ln 4})}{t^2} \,
\left[\prod_{k=5}^N \left(1 - \frac{1}{k^t}\right)\right] dt,
\label{H9}
\end{align}
and
\begin{align}
Z_N(\alpha; j) &:=
(\ln N)^{j+1} \int_0^1 N^{1-t} t^{j-2}
\left[\prod_{k=3}^N \left(1 - \frac{1}{k^t}\right)\right] dt
\nonumber
\\
&\leq (\ln N)^{j+1} \int_0^1 N^{1-t}
\left[\prod_{k=3}^N \left(1 - \frac{1}{k^t}\right)\right] dt.
\label{H10}
\end{align}
Since $R_N(\alpha; j) > 0$, from (\ref{H9}) it is clear (by bounded convergence) that
\begin{equation}
\lim_N R_N(\alpha; j) = 0.
\label{H11}
\end{equation}
To treat $Z_N(\alpha; j)$ we first use in (\ref{H10}) the estimate
$\ln(1-x) < -x$ for $0 < x < 1$, in order to deduce that
\begin{equation}
\prod_{k=3}^N \left(1 - \frac{1}{k^t}\right)
= \exp\left(\sum_{k=3}^N \ln\left(1 - \frac{1}{k^t}\right)\right)
< \exp\left(-\sum_{k=3}^N \frac{1}{k^t}\right).
\label{H12}
\end{equation}
Now,
\begin{equation*}
\sum_{k=3}^N \frac{1}{k^t} > \int_3^N \frac{d\kappa}{\kappa^t}
= \frac{N^{1-t} - 3^{1-t}}{1-t}
\end{equation*}
and hence (\ref{H12}) yields
\begin{equation}
\prod_{k=3}^N \left(1 - \frac{1}{k^t}\right)
< \exp\left(-\frac{N^{1-t} - 3^{1-t}}{1-t}\right).
\label{H13}
\end{equation}
Under (\ref{H13}), formula (\ref{H10}) implies
\begin{equation*}
Z_N(\alpha; j)
< (\ln N)^{j+1} \int_0^1 N^{1-t}
\exp\left(-\frac{N^{1-t} - 3^{1-t}}{1-t}\right) dt.
\end{equation*}
or (by substituting $s = 1-t$)
\begin{equation}
Z_N(\alpha; j)
< (\ln N)^{j+1}
\int_0^1 N^s \exp\left(-\frac{N^s - 3^s}{s}\right) ds.
\label{H14}
\end{equation}
Since we are interested in letting $N \to \infty$, we can assume that $N \geq 4$.
Then, $(N^s - 3^s) / s > (\ln N - \ln 3)$ for all $s > 0$ and hence
\begin{equation}
\int_0^{1/2} N^s \exp\left(-\frac{N^s - 3^s}{s}\right) ds
< N^{1/2} \int_0^{1/2} e^{-(\ln N - \ln 3)} ds
= \frac{3}{2} N^{-1/2}.
\label{H15}
\end{equation}
Also,
\begin{equation}
\int_{1/2}^1 N^s \exp\left(-\frac{N^s - 3^s}{s}\right) ds
< N \, \frac{1}{2} \, e^6 \exp\left(-N^{1/2} \right).
\label{H16}
\end{equation}
Therefore, since $Z_N(\alpha; j) > 0$, by using (\ref{H15}) and (\ref{H16}) in (\ref{H14}) we obtain that
\begin{equation}
\lim_N Z_N(\alpha; j) = 0.
\label{H17}
\end{equation}
Then, since $Q_N^1(\alpha; j) > 0$ by using (\ref{H17}) and (\ref{H11}) in (\ref{H8}) it follows that
\begin{equation}
\lim_N Q_N^1(\alpha; j) = 0.
\label{H18}
\end{equation}
Finally, we need to analyze $Q_N^{\infty}(\alpha; j)$ of (\ref{H3})
(recall that $a_k = \ln k$).
Since
\begin{equation*}
\frac{1}{1 - e^{-t \ln k}} \leq \frac{3}{2} < 2
\qquad \text{for all }\; k \geq 3, \; t \geq 1,
\end{equation*}
we have
\begin{equation}
\sum_{k=3}^{\infty} (\ln k)^j \frac{e^{-t \ln k}}{1 - e^{-t \ln k}}
\leq 2\sum_{k=3}^{\infty} (\ln k)^j e^{-t \ln k}
= 2\sum_{k=3}^{\infty} \frac{(\ln k)^j}{k^t}.
\label{H19}
\end{equation}
It is now easy to check that
\begin{equation}
\sum_{k=3}^{\infty} \frac{(\ln k)^j}{k^t}
\sim \int_1^{\infty} \frac{(\ln \xi)^j}{\xi^t}\, d\xi
= \frac{\Gamma(j+1)}{(t-1)^{j+1}}
= \frac{j!}{(t-1)^{j+1}}
\qquad \text{as }\; t \rightarrow 1^+.
\label{H20}
\end{equation}
Hence, from (\ref{H19}) and (\ref{H20}) it follows that there is a constant
$M > 0$ such that
\begin{equation}
\sum_{k=3}^{\infty} (\ln k)^j \frac{e^{-t \ln k}}{1 - e^{-t \ln k}}
\leq \frac{M}{(t-1)^{j+1}}
\qquad \text{for all }\; t > 1.
\label{H21}
\end{equation}
In addition, for $t \to \infty$ we have estimate (\ref{G8}), namely
\begin{equation}
\sum_{k=3}^{\infty} (\ln k)^j \frac{e^{-t \ln k}}{1 - e^{-t \ln k}}
= O\left(\frac{1}{3^t}\right).
\label{H22}
\end{equation}
Estimates (\ref{H21}) and (\ref{H22}) imply that, for $N \geq j+4$,
\begin{equation}
\int_1^{\infty} \left[\sum_{k=3}^{\infty} (\ln k)^j \,
\frac{e^{-t \ln k}}{1 - e^{-t \ln k}} \, t^{j-1} \right]
\left[\prod_{k=3}^N \left(1 - e^{-t \ln k}\right)\right] dt < \infty.
\label{H23}
\end{equation}
We, therefore, have a case similar to the one of Proposition 3: For
$N \geq j+4$ it holds that
\begin{equation*}
\left[\sum_{k=3}^N (\ln k)^j \,
\frac{e^{-t \ln k}}{1 - e^{-t \ln k}} \, t^{j-1} \right]
\prod_{k=3}^N \left(1 - e^{-t \ln k}\right)
\end{equation*}
\begin{equation}
< \left[\sum_{k=3}^{\infty} (\ln k)^j \,
\frac{e^{-t \ln k}}{1 - e^{-t \ln k}} \, t^{j-1} \right]
\prod_{k=3}^{j+4} \left(1 - e^{-t \ln k}\right)
\label{H24}
\end{equation}
and, consequently, under (\ref{H23}), dominated convergence implies
\begin{equation}
\lim_N Q_N^{\infty}(\alpha; j) =
\int_1^{\infty} \left[\sum_{k=3}^{\infty} (\ln k)^j \,
\frac{e^{-t \ln k}}{1 - e^{-t \ln k}} \, t^{j-1} \right]
\left[\prod_{k=3}^{\infty} \left(1 - e^{-t \ln k}\right)\right] dt
< \infty.
\label{H25}
\end{equation}
Finally, combining (\ref{H18}) and (\ref{H25}) with (\ref{H1}) we obtain
\begin{equation*}
\lim_N E[U_j^N] = I(\alpha; j)
= \frac{1}{(j-1)!}
\int_1^{\infty} \left[\sum_{k=3}^{\infty} \frac{(\ln k)^j}{k^t} \,
\frac{t^{j-1}}{1 - e^{-t \ln k}}\right]
\left[\prod_{k=3}^{\infty} \left(1 - \frac{1}{k^t}\right)\right] dt
< \infty.
\end{equation*}
Notice that the quantity $I(\alpha; j)$ of the above formula is the same as that of (\ref{GS4b}) (with $\alpha = \{\ln k\}_{k=3}^{\infty}$).
This completes the discussion of Example 6.\\\\

\noindent {\bf Acknowledgment.} The second author (VGP) is grateful to the Boeing Center for Technology, Information and Manufacturing (BCTIM)
of the Olin School of Business, Washington University in St. Louis, and its director, professor P. Kouvelis, for their hospitality.
This paper was written during the author's visit at BCTIM.


\end{document}